\documentclass[11pt]{article}
\usepackage{graphicx,epsfig}
\usepackage{amsmath,amssymb,euscript}
\usepackage{latexsym}
\usepackage{color}
\newtheorem{theorem}{\bf Theorem}[section]

\newcommand{\Proofend}{\hfill$\diamondsuit$}
\newcommand{\PP}{{\mathbb P}}
\newcommand{\RR}{{\mathbb R}}

\newcommand{\NN}{\mathbb{N}}

\newcommand{\C}{\mathcal{C}}
\newcommand{\V}{\mathcal{V}}
\newcommand{\W}{\mathcal{W}}
\newcommand{\bigO}{\mathcal{O}}
\newcommand{\m}{{m^{\hspace{.05cm}\mbox{\large $\bot$}}}}

\begin{document}
\title{A free parameter depending family of polynomial wavelets on a compact interval\footnote{Corresponding author: Woula Themistoclakis (woula.themistoclakis@cnr.it)}}
\author{Woula Themistoclakis\thanks{C.N.R. National
        Research Council of Italy,
        IAC Institute for Applied Computing ``Mauro Picone'',  via P. Castellino, 111, 80131 Napoli, Italy.
        woula.themistoclakis@cnr.it \\
        This author was partially supported by INdAM-GNCS}
        \and Marc Van Barel\thanks{KU Leuven, Department of Computer Science, KU Leuven,
Celestijnenlaan 200A,
B-3001 Leuven (Heverlee), Belgium. marc.vanbarel@cs.kuleuven.be.
\\
This author was partially supported by the Fund for Scientific Research–Flanders (Belgium), project G0B0123N.
}}
\maketitle

\begin{abstract}
On a compact interval, we introduce and study a whole family of wavelets depending on a free parameter that can be suitably modulated to improve performance. Such wavelets arise from de la Vall\'ee Poussin (VP) interpolation at Chebyshev nodes, generalizing previous work by Capobianco and Themistoclakis who considered a special parameter setting. In our construction, both scaling and wavelet functions are interpolating polynomials at some Chebyshev zeros of 1st kind. Contrarily to the classical approach, they are not generated by dilations and translations of a single mother function and are naturally defined on the interval $[-1,1]$ to which any other compact interval can be reduced. In the paper, we provide a non-standard multiresolution analysis with fast (DCT-based) decomposition and reconstruction algorithms. Moreover, we state several theoretical results, particularly on convergence,  laying the foundation for future applications.
\end{abstract}


\section{Introduction}
Nowadays myriads of wavelet functions are known within a well established theory that classically involves dilations and translations of a single (scaling or wavelet) mother function on the whole real axis \cite{Chui-book, Dau-book}. In adapting such a theory to a compact interval, typically various techniques are used (periodization, zero-padding etc.) that,  bringing us back to the real axis, can originate problems at the boundaries of the interval \cite{Cohen1993}. \newline
To avoid such problems, several polynomial scaling and wavelet functions have been generated without any mother to dilate and translate, but by taking particular combinations of orthogonal polynomials that naturally fit on the interval of interest. By appropriately increasing the degree of the considered basis polynomials, they give rise to a non-standard multi-resolution analysis based on the theory of orthogonal polynomials and known results in approximation theory (see, e.g., \cite{FisPre, FT-wave, KilPre, PSTasche})

In this paper we focus on particular polynomial wavelets introduced in \cite{CT-wave} starting from the de la Vall\'ee Poussin type interpolation (shortly VP interpolation) at Chebyshev zeros of the first kind \cite{Th-1999, Th-2012, Th-L1, OT-APNUM21}.

Motivated by recent applications of this type of approximation in image processing \cite{ORT-Jmiv, ORT-Etna, ORT-RemSen} as well as in the numerical quadrature of the Hilbert transform \cite{ORT-Hilbert1, ORT-Hilbert2} and in the numerical solution of integro-differential equations \cite{ORT-Cauchy, BOT-Prandtl, MT-air}, in this paper we want to present a generalization of the wavelets introduced in \cite{CT-wave}.

The main peculiarity of VP interpolation is that it depends, in addition to the number of nodes $n$, also on a free integer parameter $0<m<n$ that, appropriately chosen, can improve the pointwise approximation using the same $n$ data points (see, e.g., \cite{Th-2012, TB-GVP}).

In \cite{CT-wave} this freedom is lost by setting $m=n/2$ and consequently assuming $n$ even.
The novelty we introduce is to maintain the arbitrariness of $m$ and therefore of $n$.

In order to provide the reader with all the tools necessary for future applications, the theoretical results obtained in \cite{CT-wave}  are generalized and/or improved here, and new ones are also added. 
Finally, by generalizing the decomposition/reconstruction formulas obtained in \cite{CT-wave} , we get fast algorithms, DCT based.

The outline of the paper is the following. Section 2 collects all the notation and preliminary results we use, Section 3 is about the scaling functions and VP interpolation, Section 4 regards wavelets and their properties, Section 5 provides the decomposition and reconstruction algorithms, Section 6 contains the proofs, and finally, Section 7 deals with conclusions.
\section{Notation and preliminary results}
For any $n\in\NN$, we denote the set of the first kind Chebyshev zeros of order $n$ as follows
\[
X_n=\left\{x_k^n=\cos\frac{(2k-1)\pi}{2n},\ k=1,\ldots,n\right\}, \qquad n\in\NN .
\]
One can easily check that $X_n\subset X_{3n}$ holds $\forall n\in\NN$. Hence, setting
\[
Y_{2n}=X_{3n}-X_n=\{y_k^n, \ k=1,\ldots,2n\}, \qquad n\in\NN,
\]
we have the following decomposition
\begin{equation}\label{XY}
X_{3n}=X_n \cup Y_{2n}\qquad\mbox{and}\qquad X_n\cap Y_{2n}=\emptyset, \qquad \quad \forall n\in\NN.
\end{equation}
In the paper, we consider functions that can be sampled on the previous node sets. As functional norms, we consider 
\[ 
\|f\|_{p} =\left\{
\begin{array}{lll}
\sup_{|x|\le 1}|f(x)| & \mbox{if} & p=\infty\\ [.1in]
\displaystyle \left(\int_{-1}^1|f(x)|^p w(x)dx\right)^\frac 1p &
\mbox{if} & 1\le p<\infty
\end{array}\right.
\]
where, throughout the paper, $w$ denotes the Chebyshev weight of the first kind
\[ w(x)=\frac 1{\sqrt{1-x^2}}, \qquad |x|\le 1 .\]
Moreover, we denote by $C^0=C[-1,1]$
the space  of all continuous functions and we set 
$L^p_w=L^p_w[-1,1]=\{f \ : \ \|f\|_{p}<\infty\}$. In the case $p=2$ this is a Hilbert space equipped with the scalar product
\[
<f,g>_{L^2_w}=\int_{-1}^1 f(x)g(x)w(x)dx,
\]
Hence, in the paper, we say that the functions $f,g$ are orthogonal iff $<f,g>_{L^2_w}=0$.

Besides the previous continuous norms, for any $n\in\NN$, we also consider the discrete norms
\[
\|\vec{a}_n\|_{\ell^p}=\left\{\begin{array}{ll}
\displaystyle \left(\frac \pi{n}\sum_{k=1}^n |a_{k}|^p\right)^\frac 1p & 1\le p<\infty\\
\displaystyle \max_{1\le k\le n}|a_{k}| & p=\infty
\end{array}\right. \qquad \forall \vec{a}_n=(a_1,\ldots,a_n)\in\RR^n.
\]

In the sequel, we denote by $\PP_n$ the set of all algebraic polynomials of degree at most $n=0,1,2,\ldots$, and by
\[
E_n(f)_p=\inf_{P\in\PP_n}\|f-P\|_p,\qquad 1\le p\le\infty,\qquad n=0,1,..,
\]
the error of best (uniform or $L^p_w$) approximation of $f$ in $\PP_n$.\newline
In the special case of the uniform norm, we recall that the
Weierstrass approximation theorem ensures that
\begin{equation}\label{En-0}
\lim_{n\to\infty} E_n(f)_\infty=0, \qquad \forall f\in \C^0
\end{equation}
with a convergence order depending on the smoothness degree of $f$ \cite{DT}. In particular, we consider the class $Lip_\alpha$  of H\"older continuous functions on $[-1,1]$ with exponent $\alpha\in ]0,1]$, the class $\C^s$  of all functions that are $s$--times continuously differentiable in $[-1,1]$, and the class $\C^{s,\alpha}$ of functions $f\in \C^s$ with $f^{(s)}\in Lip_\alpha$. For these smoothness classes, we recall that, as $n\to\infty$, we have \cite{DT}
\begin{equation}\label{En-1}
E_n(f)_\infty =\left\{\begin{array}{ll}
\bigO(n^{-\alpha}) & \qquad\mbox{if}\quad f\in Lip_\alpha \quad 0<\alpha\le 1\\
\bigO(n^{-s}) & \qquad\mbox{if}\quad f\in \C^s \qquad s\in\NN\\
\bigO(n^{-s-\alpha}) & \qquad\mbox{if}\quad f\in \C^{s,\alpha} \quad\ s\in\NN, \quad 0<\alpha\le 1 .
\end{array}\right. 
\end{equation}
Throughout the paper, we denote by $p_n$ the orthonormal Chebyshev polynomial of the first kind and degree $n$, namely
\begin{equation}\label{cheb-pol}
p_n(x)=\cos(n \arccos x)\left\{\begin{array}{ll}
\sqrt{ 1/\pi} &  \mbox{ if $n=0$}\\ [.1in]
\sqrt{2/\pi} &   \mbox{ if $n>0$}
\end{array}\right. \qquad |x|\le 1, \qquad n=0,1,....
\end{equation}
and by $K_n$ the Darboux kernel
\[
K_n(x,y)=\sum_{s=0}^n p_s(x)p_s(y), \quad x,y\in [-1,1], \qquad n=0,1,....,
\]
which is the reproducing kernel in $\PP_n$ satisfying
\begin{equation}\label{repro}
\int_{-1}^1K_n(x,y)f(y)w(y)dy=f(x), \qquad\forall x\in[-1,1],\quad \forall f\in\PP_n.
\end{equation}
A crucial role in the paper is played by the following de la Vall\'ee Poussin (shortly VP) mean of the Darboux kernels \cite{Th-2012}
\begin{equation}\label{VP-mean}
v_n^m(x,y)=\frac 1{2m}\sum_{r=n-m}^{n+m-1}K_r(x,y),
\qquad n,m\in\NN,\quad m<n, \quad x,y\in [-1,1]
\end{equation}
and the related symmetric polynomial kernel
\begin{equation}\label{pol-ker-mean}
\Phi_n^m(x,y)=\frac\pi n v_n^m(x,y),
\qquad n,m\in\NN,\quad m<n, \quad x,y\in [-1,1].
\end{equation}
We recall that the VP kernel $v_n^m$ can also be interpreted as a filtered Darboux kernel of order $n+m-1$ since we have \cite{Th-2012}
\begin{equation}\label{pol-ker-sum}
v_n^m(x,y)=\sum_{r=0}^{n+m-1}\mu_{n,r}^m p_r(x)p_r(y),
\end{equation}
where
\begin{equation}\label{muj}
\mu_{n,r}^m :=\left\{\begin{array}{ll}
1 & \mbox{if}\quad 0\le r\le n-m,\\ [.1in]
\displaystyle\frac{m+n-r}{2m} & \mbox{if}\quad
n-m< r< n+m,\\
0 & \mbox{otherwise.}
\end{array}\right.
\end{equation}
Moreover, setting  $t=\arccos(x)$ and $\tau=\arccos(y)$, we have the following compact trigonometric form  \cite[p.282]{OT-APNUM21}
\begin{equation}\label{VP-trig}
v_n^m(x,y) =\frac 1{4\pi m}\left[\frac{\sin[m(t-\tau)]\sin[n(t-\tau)]}
{\sin^2[(t-\tau)/2]} + \frac{\sin[m(t+\tau)]\sin[n(t+\tau)]}
{\sin^2[(t+\tau)/2]}\right] .
\end{equation}
By using this trigonometric form, we easily deduce
\[
\int_{-1}^1|v_n^m(x,y)|w(y)dy\le \frac 1{4\pi m}\int_{-\pi}^\pi \frac{|\sin[nt]\sin[mt]|}{\sin^2[t/2]}dt \qquad |x|\le 1
\]
and applying the following estimate \cite[p.100]{Zyg-book}
\[
\frac 1{4\pi m}\int_{-\pi}^\pi \frac{|\sin[nt]\sin[mt]|}{\sin^2[t/2]}dt =\frac 4{\pi^2}\log\left[\frac{n+m-1}{2m}\right]+\bigO(1),
\]
we obtain that 
\begin{equation}\label{VP-L1}
\int_{-1}^1|v_n^m(x,y)|w(y)dy=\frac 4{\pi^2}\log\left[\frac{n+m-1}{2m}\right]+\bigO(1),\qquad \forall |x|\le 1
\end{equation}
holds for all pairs of integers $0<m<n$.
 
The VP kernel $v_n^m$  defines the VP operator 
\begin{equation}\label{VP-op}
\sigma_n^mf(x)=\int_{-1}^1v_n^m(x,y)f(y)w(y)dy,\qquad |x|\le 1,
\end{equation}
that satisfies the following invariance property
\begin{equation}\label{inva-sigma}
\sigma_n^mf=f,\qquad \forall f\in\PP_{n-m},
\end{equation}
as easily follows from  \eqref{repro}--\eqref{VP-mean}.\newline
In the sequel, we consider the Lebesgue constants of the VP operator, defined as the following operator norms
\begin{equation}\label{cp-def}
c_p=\sup_{\|f\|_p=1}\|\sigma_n^mf\|_p, \qquad 1\le p\le \infty .
\end{equation}
Note that such constants generally depend on $n$ and $m$, and the notation $c_p(n,m)$ is also used to underline this dependence when different degrees are involved (see Section \ref{Proofs}). 

In the limit cases $p=1$ and $p=\infty$, by \eqref{VP-L1}, we have
\begin{equation}\label{c_inf}
c_1=c_\infty=\max_{|x|\le 1}\int_{-1}^1|v_n^m(x,y)|w(y)dy=\frac 4{\pi^2}\log\left[\frac{n+m-1}{2m}\right]+\bigO(1),
\end{equation}
while the interpolation theorem in weighted $L^p$ spaces (cf. \cite[Corollary 2.2]{Gus}, see also \cite{TB-GVP}) yields
\begin{equation}\label{LCp}
c_p\le const. \ \cdot c_\infty, \qquad 1<p<\infty.
\end{equation}

\section{Scaling functions and VP interpolation}
Contrary to the standard definition, we do not take a unique function to dilate and translate for generating the scaling functions,  but for arbitrarily fixed $n,m\in\NN$ with $n>m$,  we set the following $n$ polynomials
\begin{equation}\label{sca}
\Phi_{n,k}^m(x)=\Phi_n^m(x_k^n,\ x)=\frac\pi n\ v_n^m(x_k^n, \ x),\qquad k=1,\ldots,n\qquad x\in[-1,1].
\end{equation}
as {\it scaling functions of order $n$ and parameter $m$}.

We recall that these polynomials have already been introduced in the literature as {\it fundamental VP polynomials} \cite{Th-2012}, and have been studied in several papers (see \cite{OT-APNUM21, OT-DRNA21} and the references therein). \newline
Similarly to fundamental Lagrange polynomials interpolating at $X_n$, they are linearly independent and satisfy the interpolation property \cite[Eq. (2.10)]{CT-wave}
\begin{equation}\label{sca-int}
\Phi_{n,k}^m(x_h^n) =\delta_{h,k}=\left\{\begin{array}{ll}
1 & h=k\\
0 & h\ne k
\end{array}\right.\qquad h,k=1,\ldots,n,\qquad 0<m<n.
\end{equation}
Besides $k$, they depend on the integers $0<m<n$ where $n$ indicates the resolution degree, i.e., the number of the interpolation nodes in $X_n$, and $m$ is a free parameter that can be arbitrarily chosen in $\{1,2,\ldots, n-1\}$ without destroying the interpolation property \eqref{sca-int}. 

\begin{figure}[!htb]%
\begin{center}
\includegraphics[scale = 0.45]{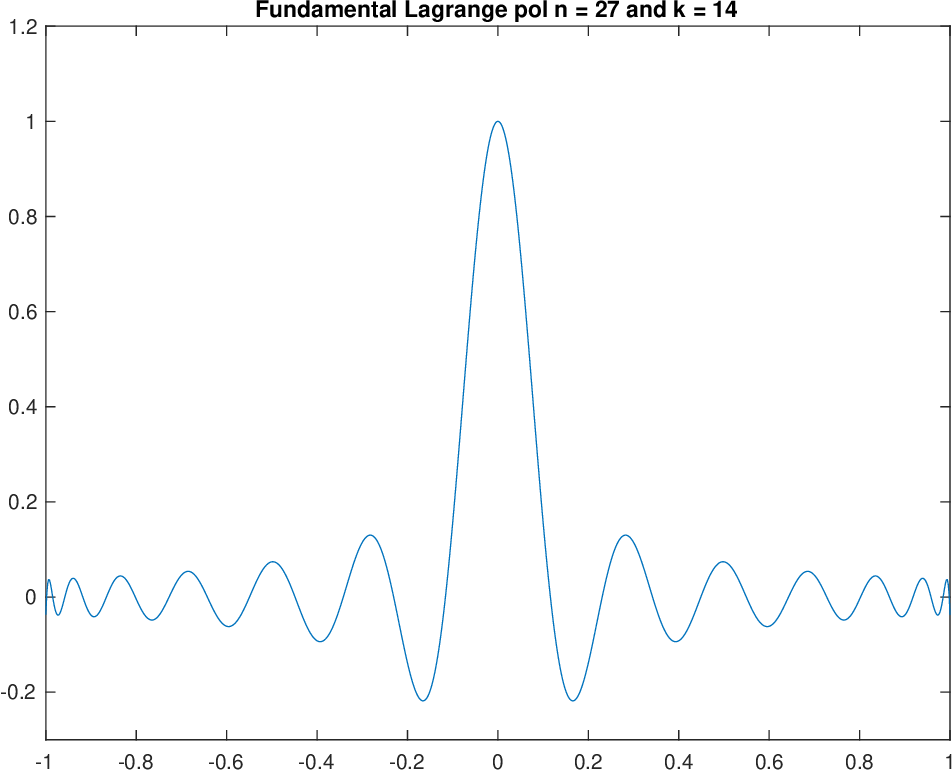}
\includegraphics[scale=0.45]{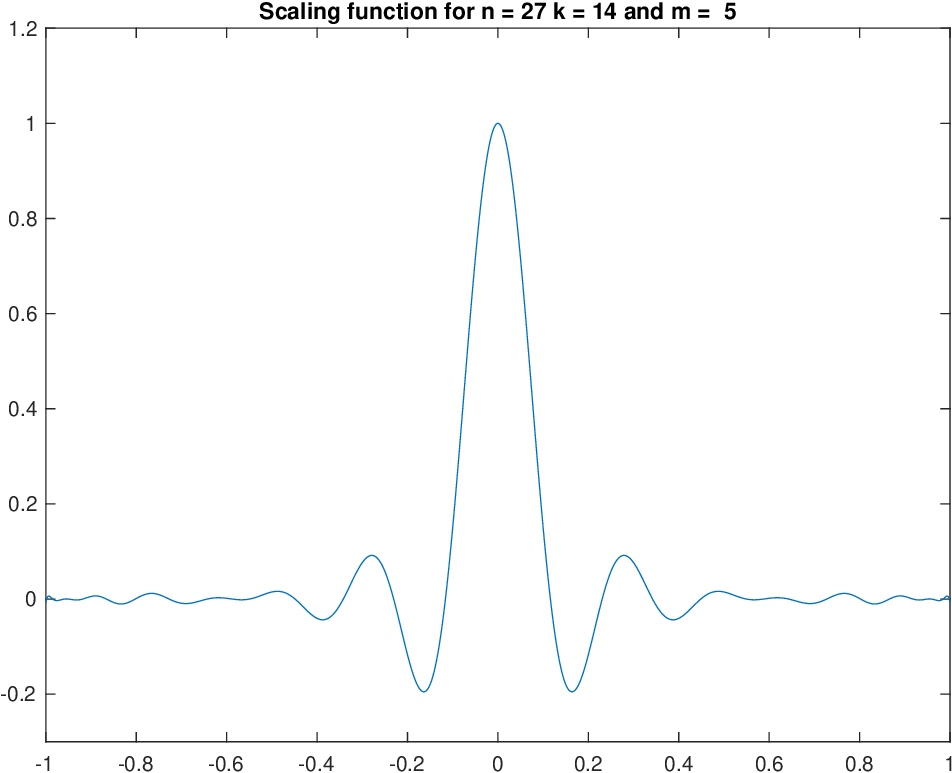}\\
\includegraphics[scale=0.45]{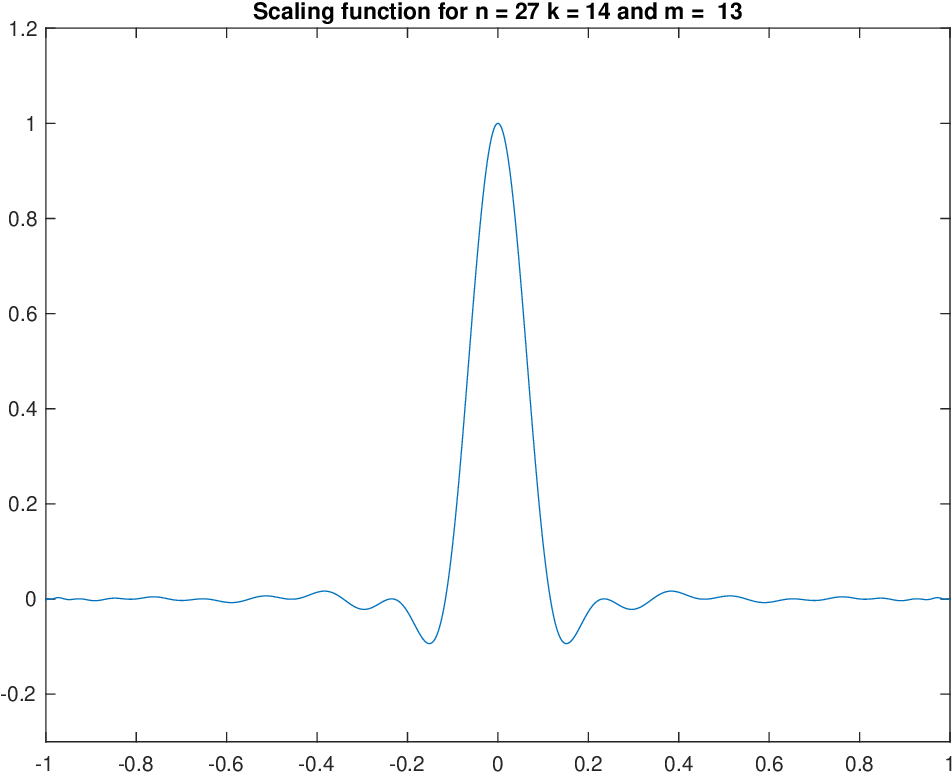}
\includegraphics[scale=0.45]{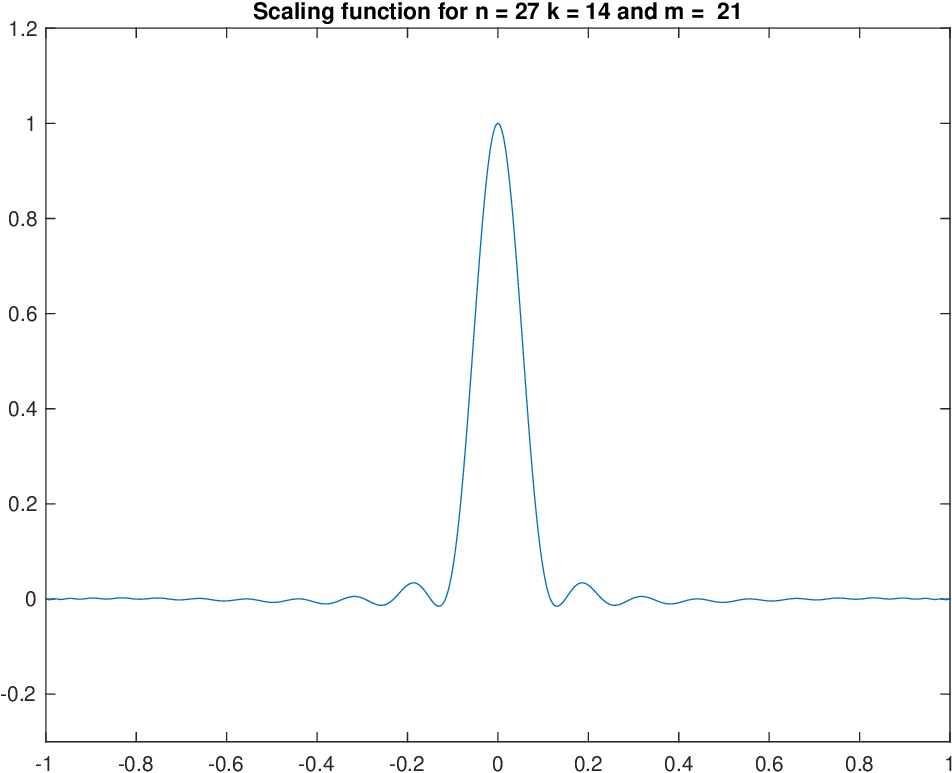}
\end{center}
\caption{ The oscillating behavior of the fundamental Lagrange polynomial compared to the fundamental VP polynomials with several $m$.\label{fig_compare_Lagr_VP} }
\end{figure}

\begin{figure}[!htb]%
\begin{center}
\includegraphics[scale = 0.45]{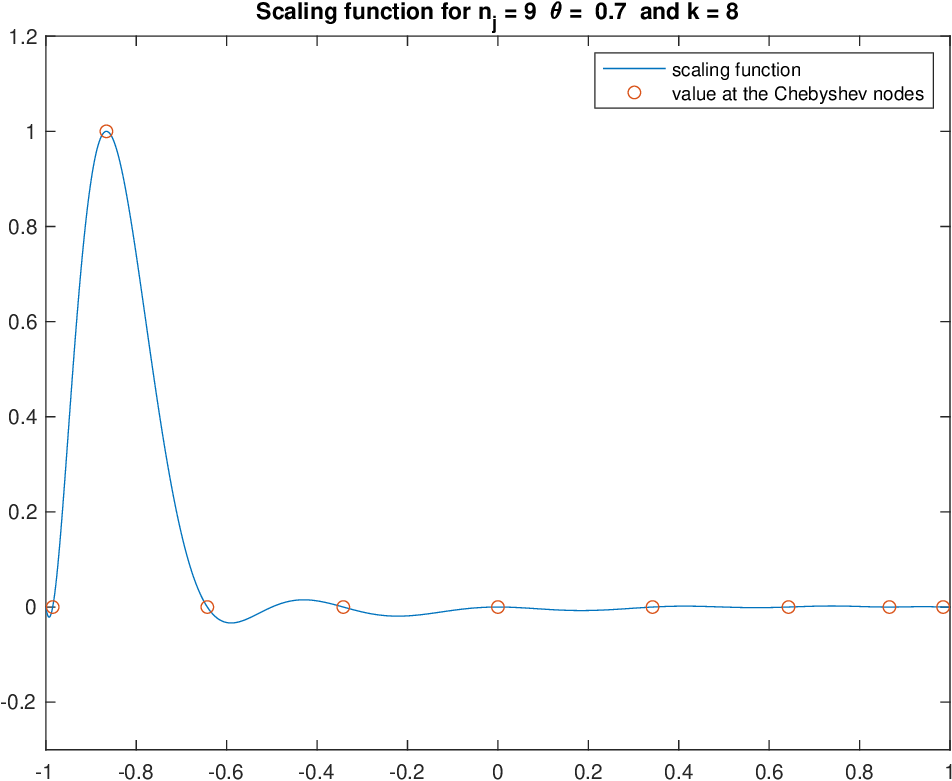}
\includegraphics[scale = 0.45]{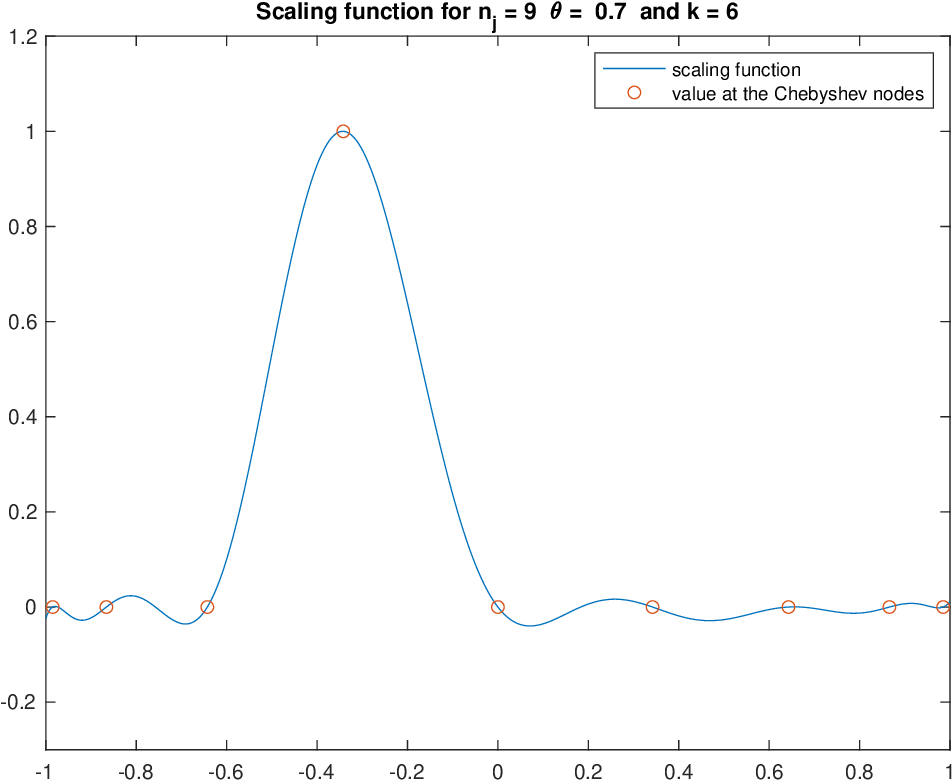}\\
\includegraphics[scale = 0.45]{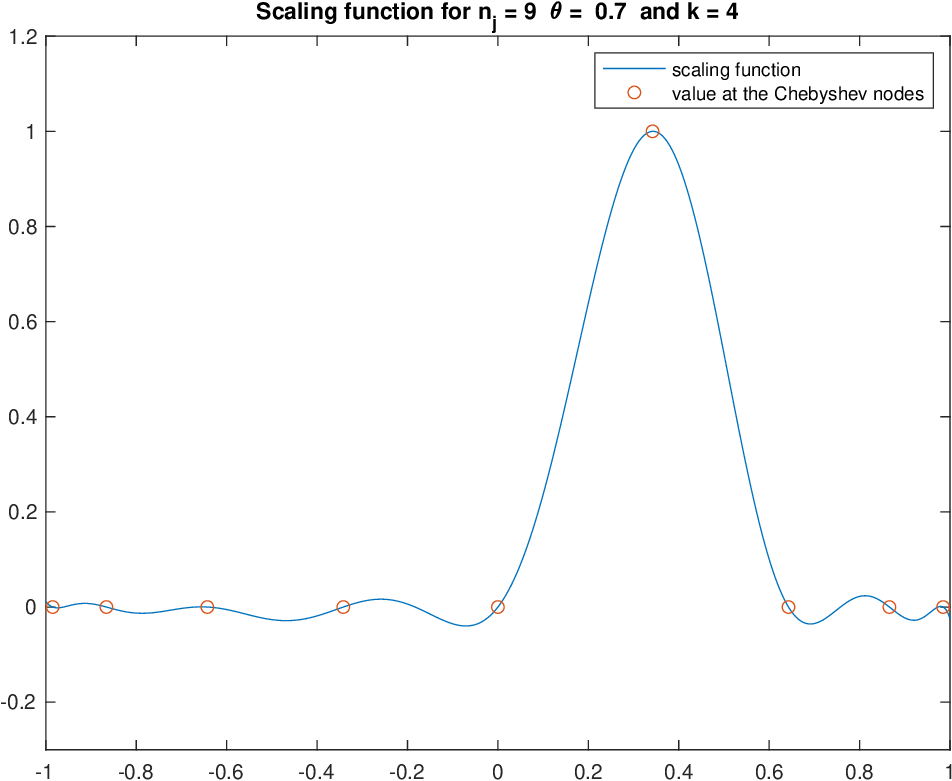}
\includegraphics[scale = 0.45]{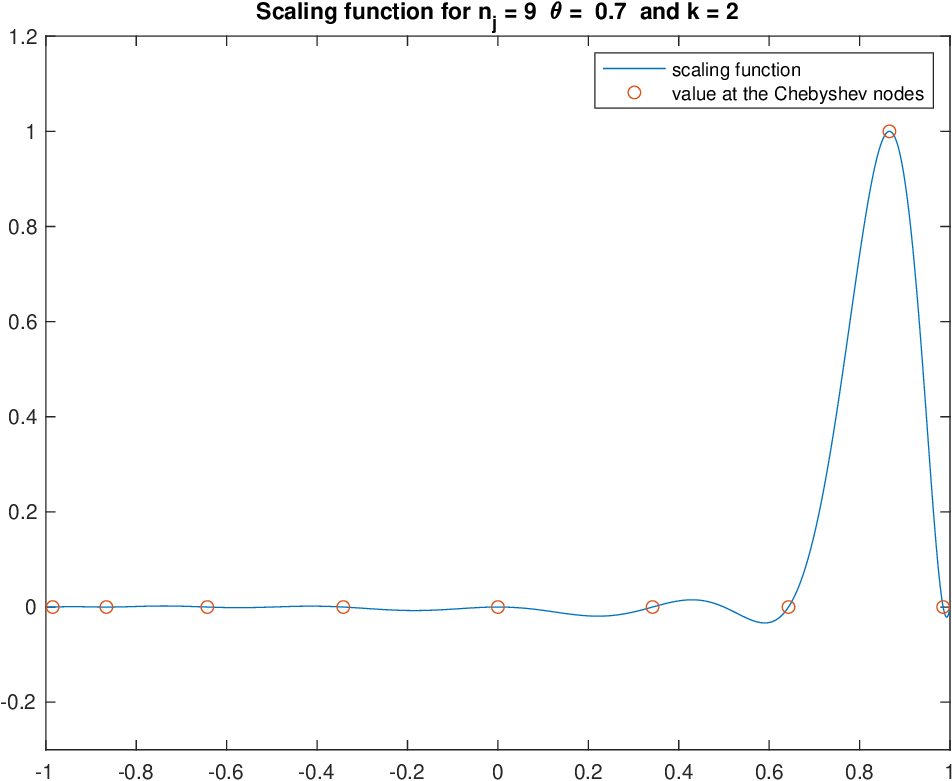}
\end{center}
\caption{ The scaling function in function of $k$ with $m=\lfloor \theta n\rfloor$.\label{fig_scaling_functions_k} }
\end{figure}

\begin{figure}[!htb]%
\begin{center}
\includegraphics[scale = 0.45]{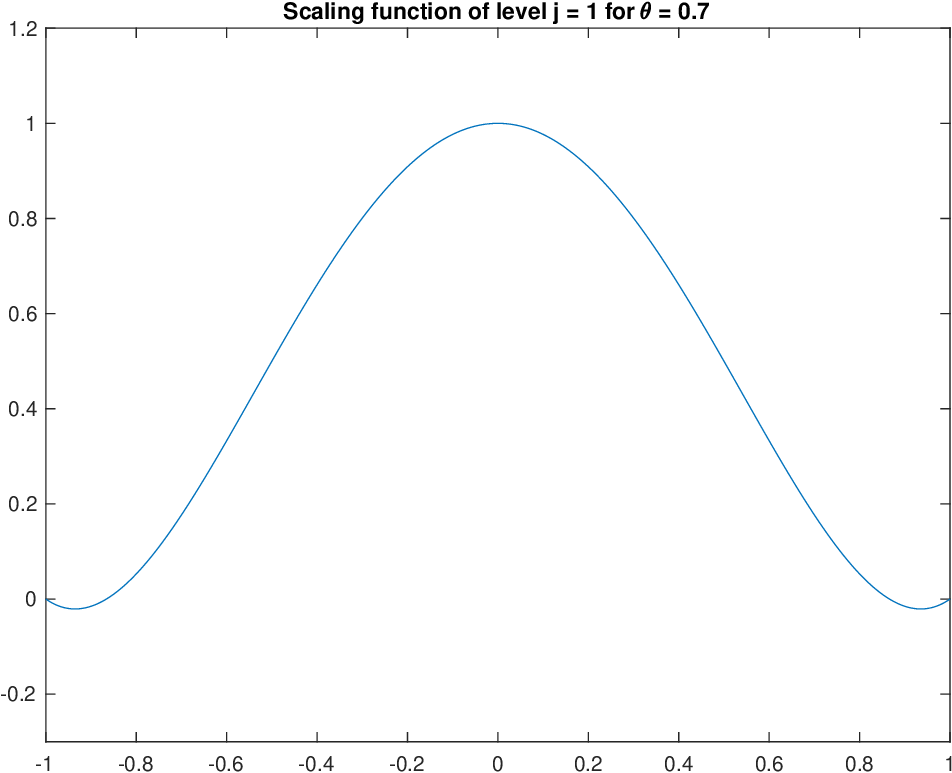}
\includegraphics[scale = 0.45]{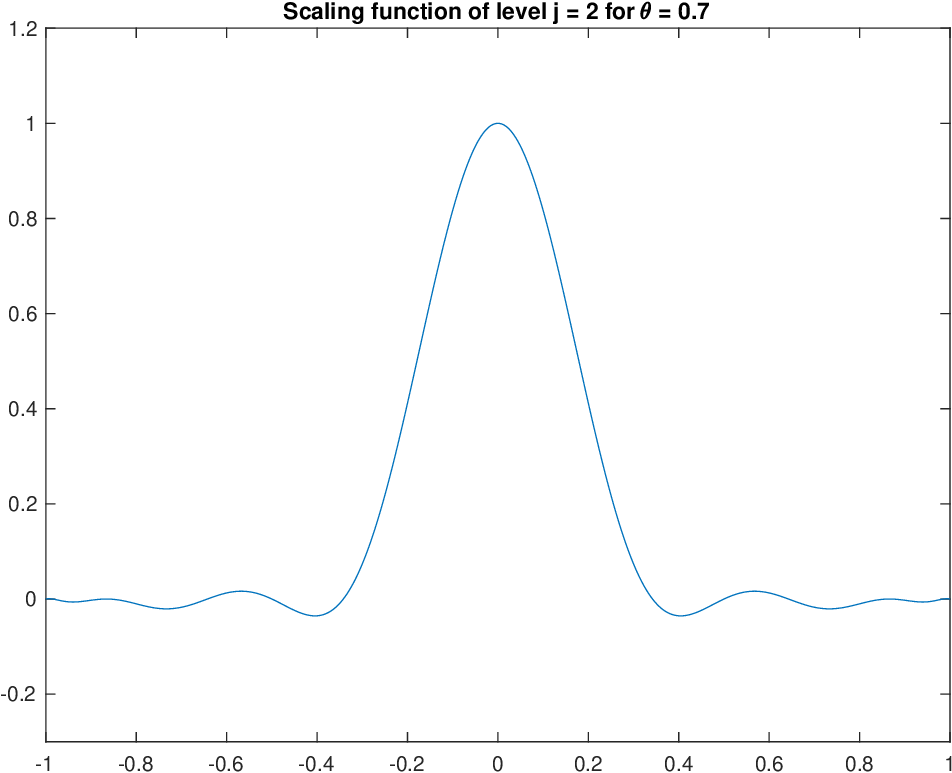}\\
\includegraphics[scale = 0.45]{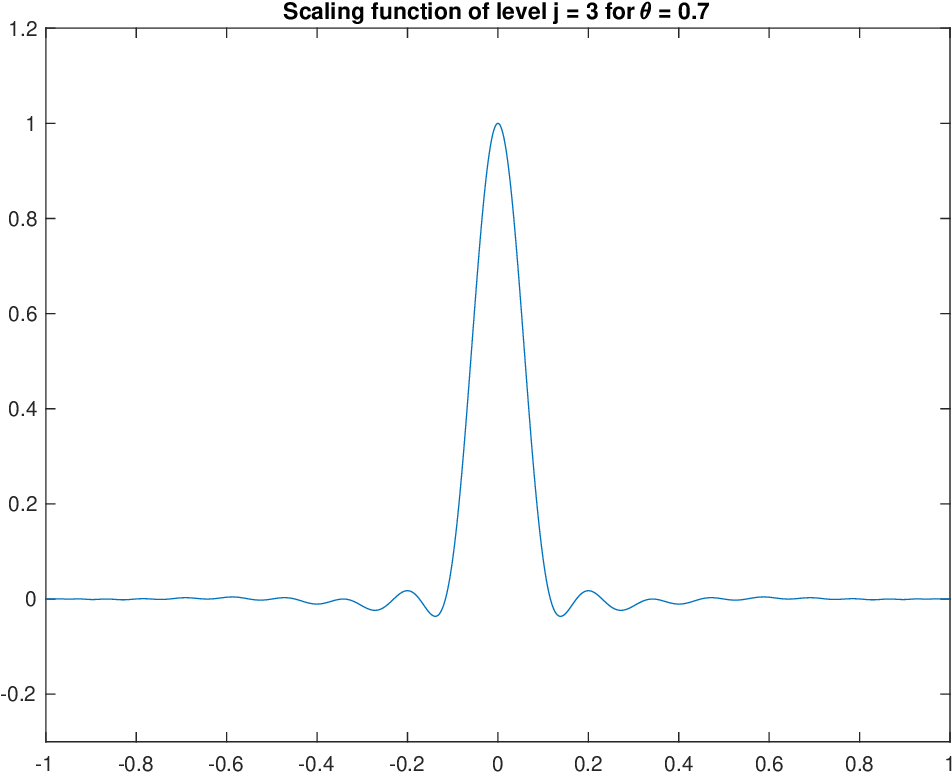}
\includegraphics[scale = 0.45]{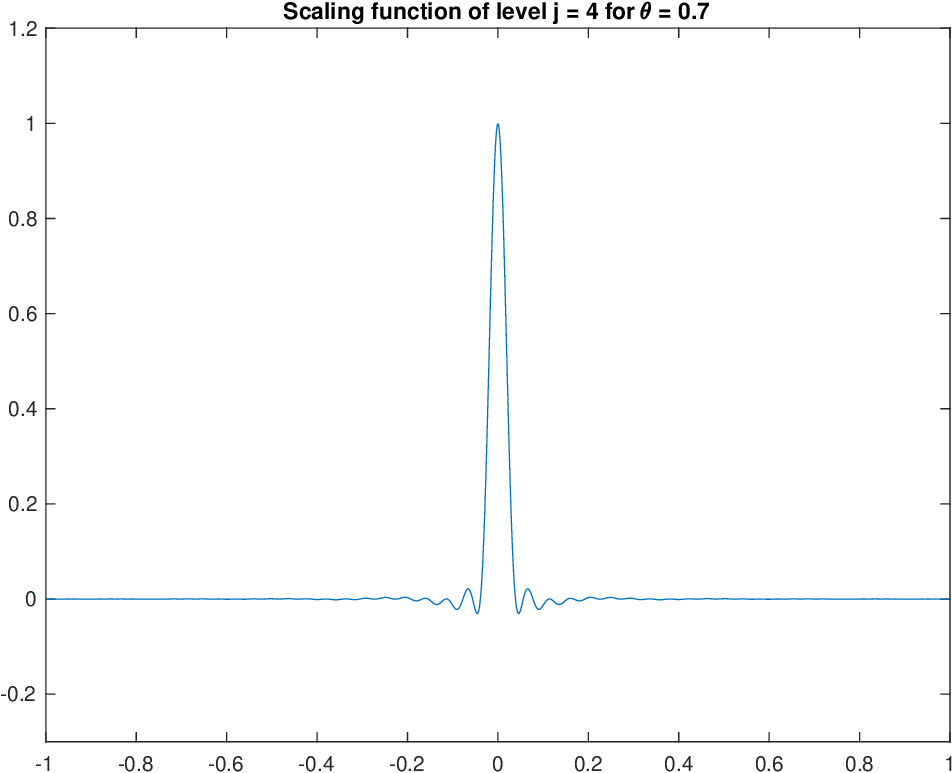}
\end{center}
\caption{ Scaling functions for different values of $n = 3^j$ with $m=\lfloor \theta n\rfloor$.\label{fig_scaling_function_n} }
\end{figure}

In Figure \ref{fig_compare_Lagr_VP} we see how this additional parameter can soften the typical oscillating behavior of the fundamental Lagrange polynomial.
Moreover, by Figure \ref{fig_scaling_functions_k} , we see each scaling function $\Phi_{n,k}^m$ is well localized around the Chebyshev node $x_k^n$ where it assumes its maximum value (=1). 
Finally, Figure \ref{fig_scaling_function_n} shows how such localization increases as $n$ increases. 

Other properties of the scaling functions are summarized in the following theorem, where we point out that $(e)$ improves a recent estimate stated in \cite[Thm. 6.1]{OT-DRNA21}.
\begin{theorem}\label{th-sca}
For all $n,m\in\NN$ with $n>m$, we have the following
\begin{itemize}
\item[(a)] $\displaystyle \int_{-1}^1 x^s\Phi_{n,k}^m(x)w(x)dx=\frac\pi n (x_k^n)^s,\qquad s=0,1,\ldots, n-m,\qquad k=1,\ldots,n$\\
\item[(b)] $\displaystyle x^s=\sum_{k=1}^n (x_k^n)^s\Phi_{n,k}^m(x),\qquad  s=0,1,\ldots, n-m, \quad  x\in [-1,1]$\\
\item[(c)] $\displaystyle \|\Phi_{n,k}^m\|_\infty =1, \qquad\quad k=1,\ldots,n$\\
\item[(d)] $\displaystyle  \frac\pi{n^p}\le\ \|\Phi_{n,k}^m\|_p^p\ \le \frac\pi n\ c_\infty,\qquad 1\le p<\infty,\qquad\quad k=1,\ldots,n$\\
\item[(e)]  $\displaystyle \sum_{k=1}^n|\Phi_{n,k}^m(x)|\le\left[1+2\pi\ \frac{n+m-1}{n}\right]\ c_\infty\le (1+4\pi) c_\infty \qquad |x|\le 1$\\
\item[(f)]  For any $\vec{a}_n=(a_{1}, \ldots, a_{n})\in\RR^n$,  $1\le p<\infty$, and $\frac 1p+\frac 1{p'}=1$ we have
\begin{eqnarray}
\label{Riesz-inf}
\hspace{-1cm} \|\vec{a}_n\|_{\ell^\infty}\le & \left\|\sum_{k=1}^{n} a_{k}\Phi_{n,k}^m\right\|_\infty &\le 
c_\infty\left[1+2\pi\ \frac{n+m-1}{n}\right] \|\vec{a}_n\|_{\ell^\infty}\\
\label{Riesz-p}
\hspace{-1cm}\frac 1{1+4\pi }\|\vec{a}_n\|_{\ell^p}\le &\left\|\sum_{k=1}^{n} a_{k}\Phi_{n,k}^m\right\|_p&\le c_{p'} (1+4\pi )\ \|\vec{a}_n\|_{\ell^p}.
\end{eqnarray}
\end{itemize}
\end{theorem}
The scaling functions of a given order $n\in\NN$ and parameter $m<n$ generate the following approximation space
\begin{equation}
\label{sca-V}
\V_n^m = span \{\phi_{n,k}^m, \ k=1,\ldots,n\}, \qquad  \dim \V_n^m=n.
\end{equation}
By the previous properties, this space is nested between classical polynomial spaces as follows
\begin{equation}\label{nest}
\PP_{n-m}\subset \V_n^m\subset \PP_{n+m-1}, \qquad \forall n,m\in\NN, \quad m<n.
\end{equation}
Moreover, the scaling functions in \eqref{sca-V} constitute an interpolating Riesz basis of $\V_n^m$. 

It is also known that an orthogonal basis of $\V_n^m$ is given by the following polynomials \cite[Thm. 2.2]{CT-wave}
\begin{equation}\label{q-basis}
\Phi_{n,r}^{m^{\hspace{.05cm}\mbox{\large $\bot$}}}(x):=\left\{\begin{array}{ll}
p_r(x)&\mbox{if}\quad 0\le r\le
n-m,\\ [.1in]
\displaystyle\mu_{n,r}^m p_r(x)-\mu_{n,2n-r}^m p_{2n-r}(x)
&\mbox{if}\quad n-m<r<n.
\end{array}\right.
\end{equation}
Setting 
\begin{equation}\label{q-prod}
\nu_{n,r}^m=<\Phi_{n,r}^\m,\ \Phi_{n,r}^\m>_{L^2_w}=\left\{\begin{array}{ll} 1 & \mbox{if}\quad 0\le r\le
n-m,\\ [.1in] \displaystyle\frac{m^2+(n-r)^2}{2m^2} &
\mbox{if}\quad
 n-m<r< n,
\end{array}\right.
\end{equation}
we have
\[
<\Phi_{n,r}^\m,\ \Phi_{n,s}^\m>_{L^2_w}=\delta_{r,s} \cdot \nu_{n,r}^m,  \qquad r,s=0,..., n-1.
\]
Moreover,  the change of bases is given by the following formula \cite[Eq. (2.15)]{CT-wave}
\begin{equation}\label{qbasis-trans}
\Phi_{n,k}^m(x)=\frac\pi n\sum_{r=0}^{n-1}p_r(x_k^n)\Phi_{n,r}^{m^{\hspace{.05cm}\mbox{\large $\bot$}}}(x),
\qquad k=1,\ldots,n.
\end{equation}
Due to \eqref{sca-int}, any function $f$, sampled at $X_n$, can be projected onto the approximation space $\V_n^m$ by taking the polynomial
\begin{equation}\label{Vnm}
V_n^m f(x)=\sum_{k=1}^nf(x_k^n)\phi_{n,k}^m(x),\qquad n>m\in\NN,\qquad x\in [-1,1]
\end{equation}
that interpolates $f$ at $X_n$, i.e. 
\begin{equation}\label{interp-V}
V_n^mf(x_h^n)=f(x_h^n), \qquad h=1,\ldots,n,\qquad n>m\in\NN.
\end{equation}
Also, note that we trivially have
\begin{equation}\label{V-equi}
f\in \V_n^m \Longleftrightarrow  f=V_n^m f
\end{equation}
and in particular
\begin{equation}\label{inva}
V_n^m P=P, \qquad \forall P\in\PP_{n-m}.
\end{equation}
Moreover, by \eqref{qbasis-trans},  the polynomial  \eqref{Vnm} can be expanded in terms of the orthogonal functions \eqref{q-basis} as follows
\begin{equation}\label{Vnm-q}
V_n^m f(x)=\sum_{r=0}^{n-1}c_{n,r}\Phi_{n,r}^{m^{\hspace{.05cm}\mbox{\large $\bot$}}}(x), \qquad c_{n,r}=\left[\frac\pi n\sum_{k=1}^nf(x_k^n)p_r(x_k^n)\right], \qquad \forall n>m\in\NN .
\end{equation}
In the literature, this polynomial is already known as the VP polynomial interpolating $f$ and its approximation properties have been studied in several papers (see  \cite{OT-APNUM21, OT-DRNA21} and the references therein). 

Here, we recall that it offers several advantages over the classical Lagrange polynomial interpolating $f$ at the same nodes, especially in the case of functions a.e. smooth with isolated singularities. In such a case, one can take advantage of the presence of the free parameter $m<n$ that can be suitably modulated to attenuate the Gibbs phenomenon and improve the point-wise approximation keeping the number of nodes fixed \cite{Th-2012}.

As regards the uniform approximation, from \eqref{inva} and Thm.~\ref{th-sca} we easily deduce the following estimate of the theoretical error
\begin{equation}\label{err-VP}
E_{n+m-1}(f)_\infty\le \|f-V_n^mf\|_\infty\le \C E_{n-m}(f)_\infty,
\end{equation}
where $\C:=1+(1+4\pi)c_\infty$.

Moreover, if we are using corrupted data, say 
\[
\tilde f(x_k^n):=f(x_k^n)+\epsilon_k^n,\qquad k=1,\ldots,n
\]
then, by Thm.\ref{th-sca}, we get
\begin{equation}\label{cond-VP}
\frac 1{\C}\|\vec{\epsilon}_n\|_{\ell^p}\le \|V_n^m f-V_n^m\tilde f\|_p\le \C\|\vec{\epsilon}_n\|_{\ell^p},\qquad 1\le p\le \infty,
\end{equation}
where $\C:=(1+4\pi)c_{p'}$, $\frac 1p+\frac1{p'}=1$, and $\vec{\epsilon}_n=(\epsilon_1^n,\ldots,\epsilon_n^n)$.

In the previous estimates \eqref{err-VP}--\eqref{cond-VP}, we remark that, contrary to $n$ which is determined by the number of nodes, the parameter $m$ can be freely chosen in $\{1,\ldots,n-1\}$. 
But, as $n\to\infty$, the asymptotic estimates \eqref{c_inf}--\eqref{LCp} imply a logarithmic growth of the constant $\C$ in \eqref{err-VP} --\eqref{cond-VP}, unless we take $m$ growing at the same order. Indeed, we have
\begin{equation}\label{LCnm}
1< \frac nm=\bigO(1) \ \Longrightarrow \ c_p=\bigO(1),\quad 1\le p\le \infty,\qquad \quad \forall n,m\in\NN.
\end{equation} 

In the literature, a common choice is to take $m=\lfloor \theta n\rfloor$ with $\theta\in ]0,1[$ arbitrarily fixed. In such a case, in \eqref{err-VP} --\eqref{cond-VP}, $\C$ is an absolute constant (only depending on $\theta$) and we also have that $(n\pm m)$ behaves as $n$. Hence, by \eqref{err-VP} and \eqref{En-0},  we get
\begin{equation}\label{lim-wav}
\lim_{\scriptsize \begin{array}{c}
n\to\infty\\
m=\lfloor \theta n\rfloor\end{array}}
\|f-V_n^mf\|_\infty=0, \qquad \forall f\in C^0,\qquad \forall\theta\in ]0,1[,
\end{equation}
where the convergence order is comparable with the convergence rate of $E_n(f)_\infty$.

In particular, recalling \eqref{En-1}, for arbitrarily large $n\in\NN$ and $m=\lfloor \theta n\rfloor$ with  fixed $\theta\in ]0,1[$, we have
\begin{equation}\label{err-Vnm}
\|f-V_n^mf\|_\infty\le\C  \left\{
\begin{array}{ll}
 n^{-\alpha} & \mbox{if}\quad  f\in Lip_\alpha \qquad 0<\alpha\le 1 \\
 n^{-s}& \mbox{if}\quad  f\in C^s \qquad\quad s\in\NN\\
 n^{-s-\alpha}& \mbox{if}\quad  f\in C^{s,\alpha} \qquad s\in\NN, \quad 0<\alpha\le 1 
 \end{array}\right.
\end{equation}
where $\C>0$ is a constant independent of $n,m$ but depending on $\theta\in ]0,1[$ and $f$.

Finally, we recall that by using the samples of $f$ at $X_n$ we can simultaneously approximate also the derivatives of $f$ in $[-1,1]$ by means of the derivatives of the VP polynomial $V_n^mf$ given by
\[
(V_n^mf)^{(r)}(x)=\sum_{k=1}^n f(x_k^k) \left(\Phi_{n,k}^m\right)^{(r)}(x), \qquad |x|\le 1.
\]
For brevity, we omit the details about the computation of the derivatives of the scaling function basis, referring the interested reader to \cite{OT-DRNA21}. Here, we only recall the error estimate proved in \cite[Thm. 5.3]{OT-DRNA21}.
\begin{theorem}\label{th-Vnm}
Let $n,m\in\NN$ be arbitrarily large but such that $n/m\in [a, b]$ holds, where $b\ge a>1$ are independent of $n$ and $m$. If $f\in C^{2s}$ for some 
$s\in\NN$ then we have   
\begin{equation}\label{V-der}
 \displaystyle \|f^{(r)}-(V_n^mf)^{(r)}\|_\infty\le \frac \C{n^{2s-r}}, \qquad r=0,1,\ldots,s
\end{equation}
where $\C>0$ is independent of $n,m$.
\end{theorem}


\section{Wavelet functions and detail spaces}
The starting point for the construction of our wavelets is given by the following inclusion
\begin{equation}\label{incl}
\V^m_{n}\subset \V^m_{3n}, \qquad \forall n>m\in\NN,
\end{equation}
which easily follows from \eqref{nest} and $\PP_{n+m-1}\subset \PP_{3n-m}$.

By taking into account \eqref{incl}, for all integers $n,m\in\NN$ with $n>m$ we define the {\it detail or wavelet space} $\W_n^m$ as the orthogonal complement of $\V^m_n$ in $\V_{3n}^m$, namely
\begin{equation}\label{v+w}
\V^m_ {3n}=\V^m_{n}\oplus \W_n^m \quad\mbox{and}\quad \V_n^m\bot\W_n^m, \qquad\quad \forall n>m\in\NN.
\end{equation}
Wavelet functions are generally defined as localized basis functions that generate the detail spaces. Since $\dim \V_n^m=n$ implies that $\dim \W_n^m=2n$, we expect $2n$ wavelets for any resolution degree $n\in\NN$ and any integer parameter $0<m<n$.

Recalling \eqref{XY} and \eqref{sca-int}, it comes natural to define the {\it interpolating wavelet functions of order $n$ and parameter $m$} as those basis functions of $\W_n^m$, say
\begin{equation}\label{wav-W}
\W_n^m=span \{\psi_{n,k}^m, \ k=1,\ldots,2n\},
\end{equation}
satisfying the following interpolation property on the nodes set $Y_{2n}=X_{3n}-X_n$
\begin{equation}\label{wav-int}
\psi_{n,k}^m(y_h^n) =\delta_{h,k},\qquad h,k=1,\ldots,2n,\qquad 0<m<n.
\end{equation}
We note that \eqref{wav-W} and \eqref{wav-int} uniquely determine the interpolating wavelets $\psi_{n,k}^m$. Their explicit form is given below in terms of the polynomial kernel introduced in \eqref{pol-ker-mean}.
\begin{theorem}\label{th-wav}
For all $n,m\in\NN$ with $n>m$, the interpolating wavelets defined by \eqref{wav-W} and \eqref{wav-int} have the following analytic form
\begin{equation}\label{wav}
\psi_{n,k}^m(x)=\Phi_{3n}^m(y_k^n,x) - \sum_{h=1}^n \Phi_{n,h}^m(y_k^{n})\Phi_{3n}^m(x_h^n,x), \qquad k=1,\ldots,2n.
\end{equation}
\end{theorem}
We point out that similarly to the scaling functions, the wavelets are also not generated from translations and dilations of a mother wavelet. The notation $\psi_{n,k}^m$ underlines their dependence on the resolution level $n$, on the free parameter $m$, and on the node $y_k^n$ around which they are maximally located. 

\begin{figure}[!htb]%
\begin{center}
\includegraphics[scale = 0.45]{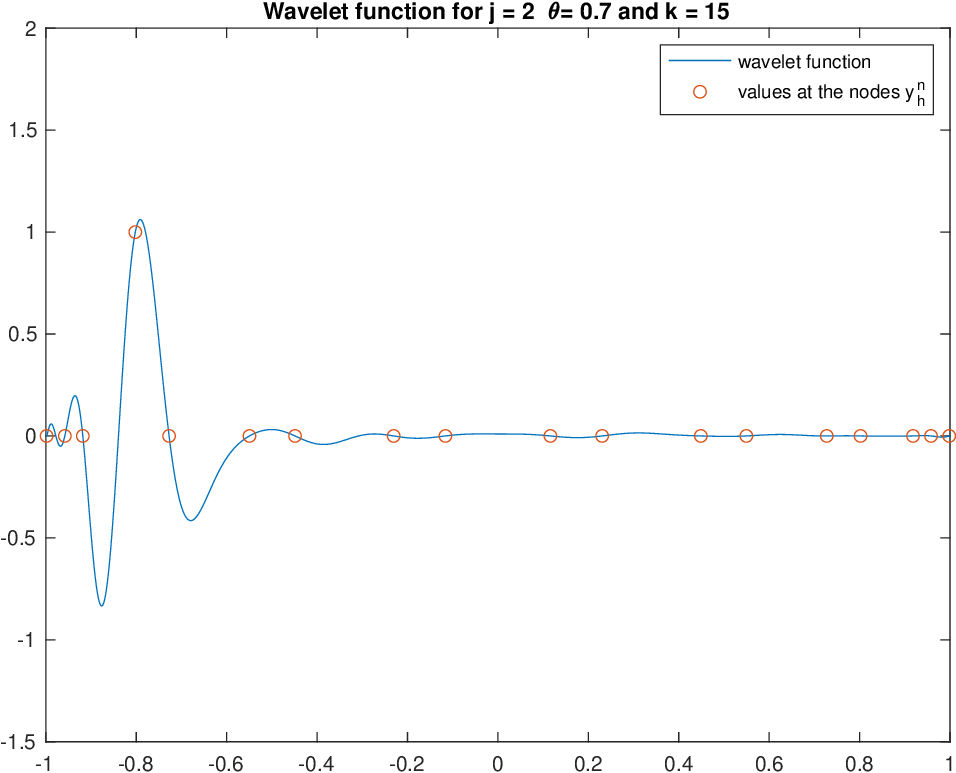}
\includegraphics[scale = 0.45]{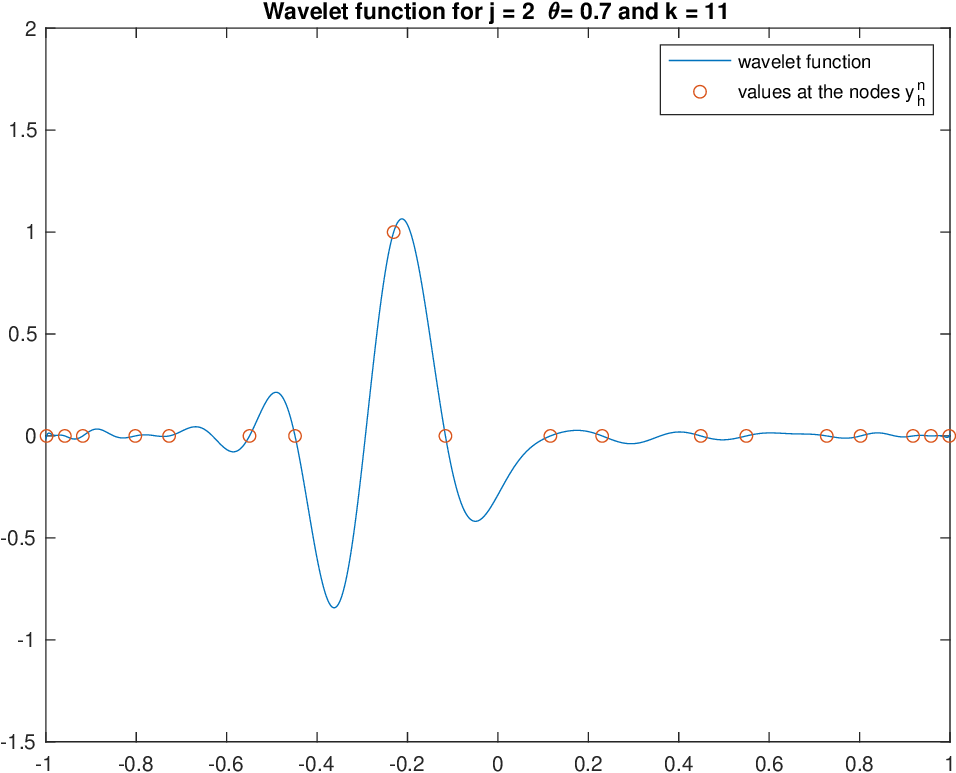}\\
\includegraphics[scale = 0.45]{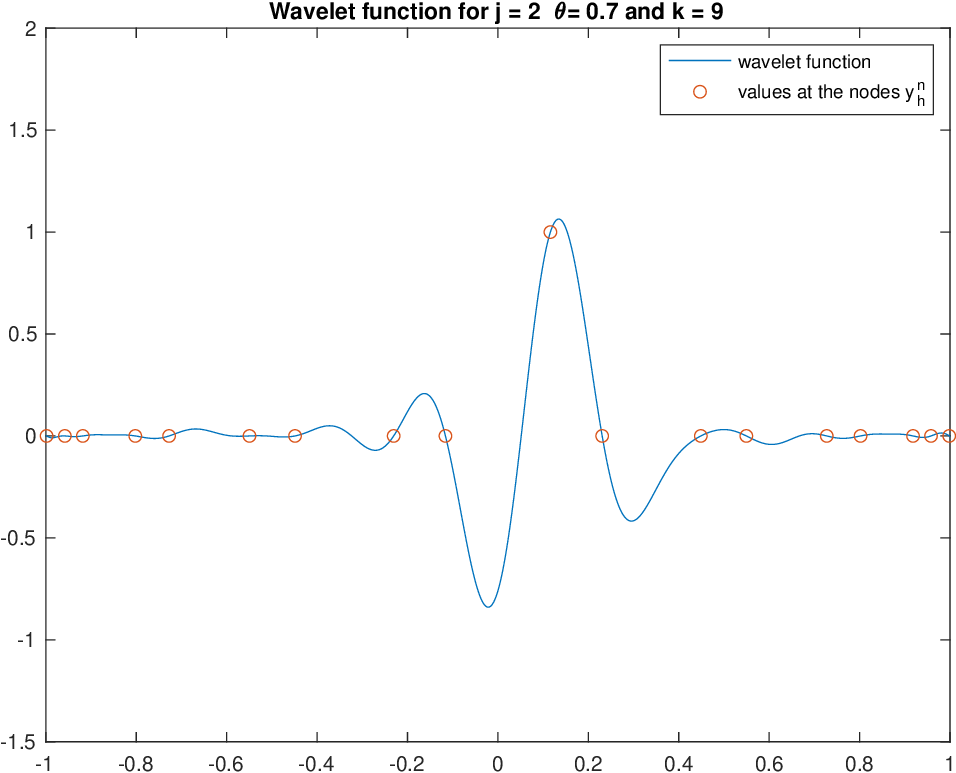}
\includegraphics[scale = 0.45]{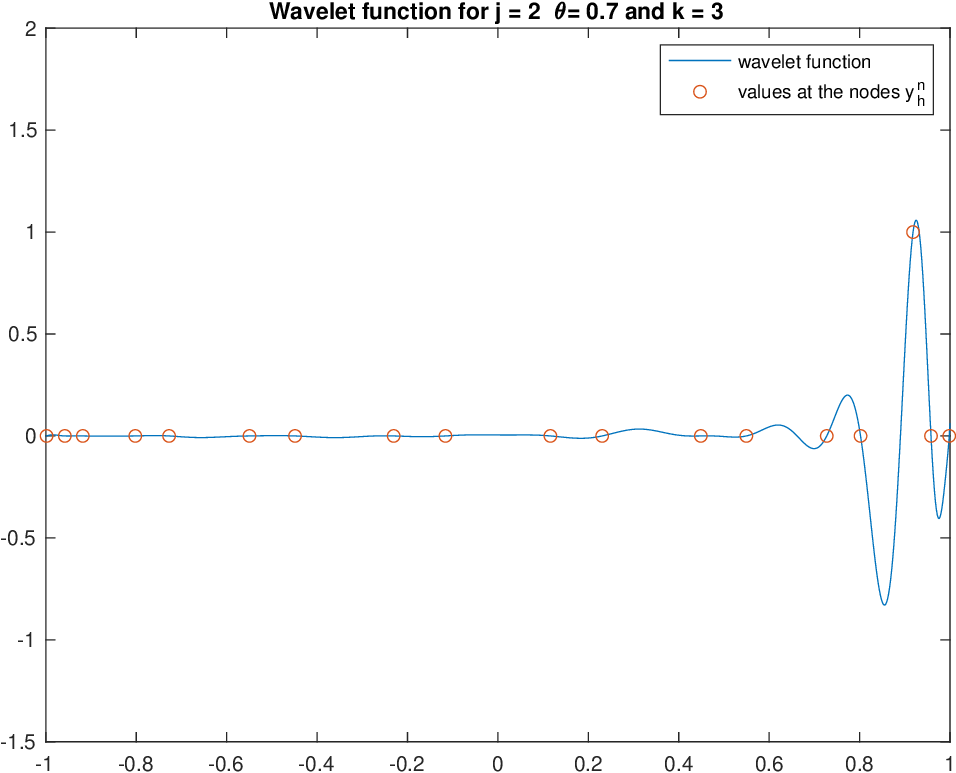}
\end{center}
\caption{ The interpolating wavelet function in function of $k$.\label{fig_wavelet_functions_k} }
\end{figure}

\begin{figure}[!htb]%
\begin{center}
\includegraphics[scale = 0.45]{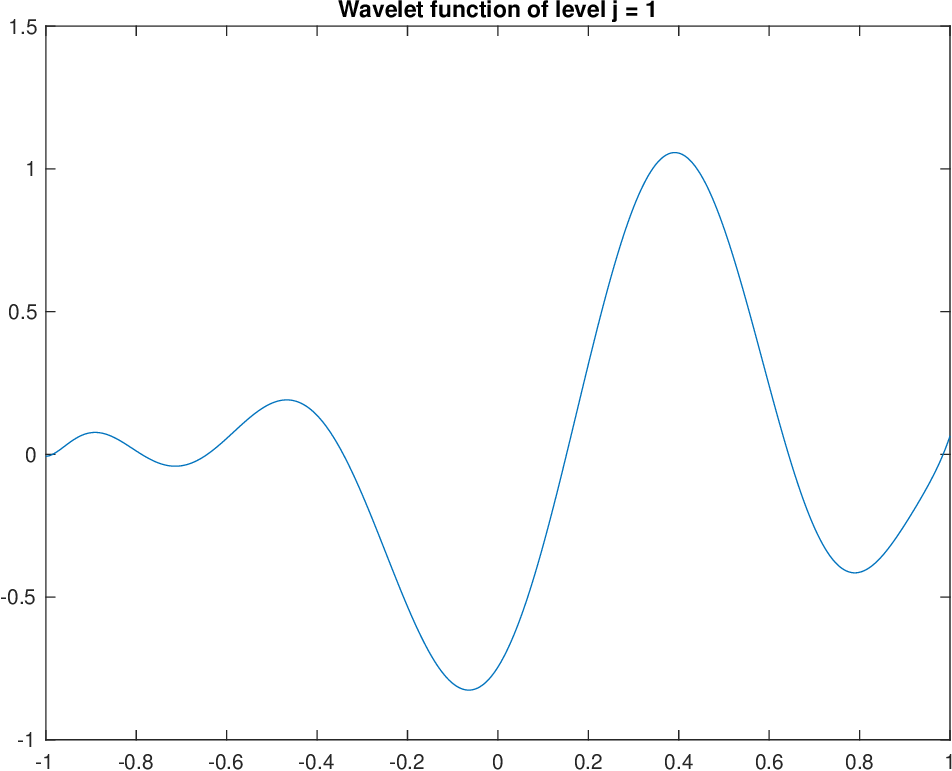}
\includegraphics[scale = 0.45]{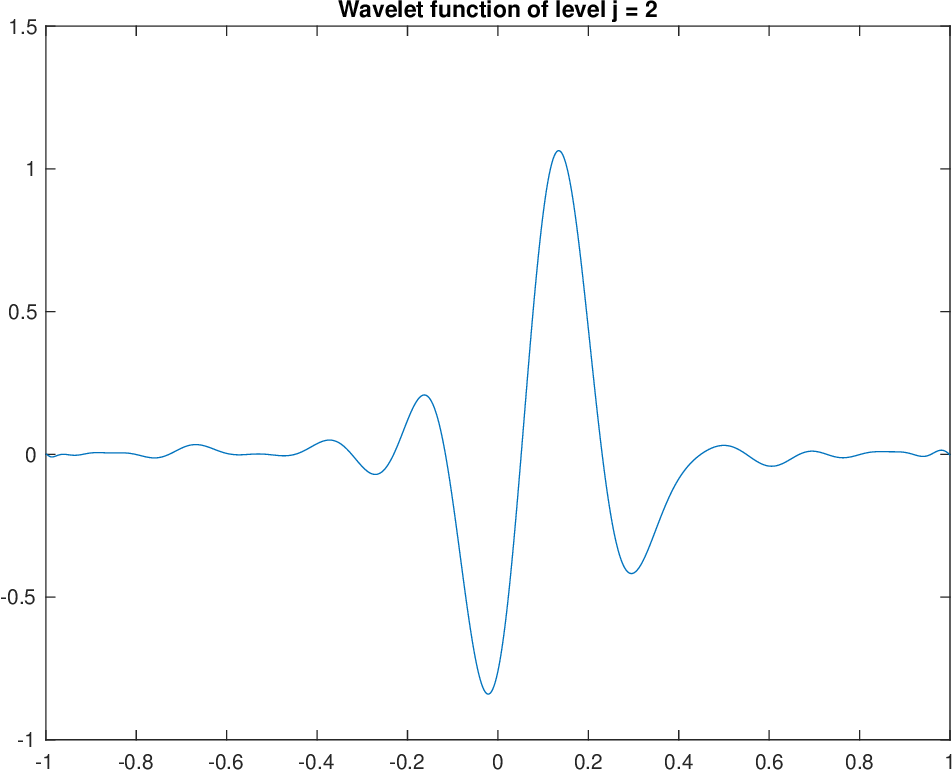}\\
\includegraphics[scale = 0.45]{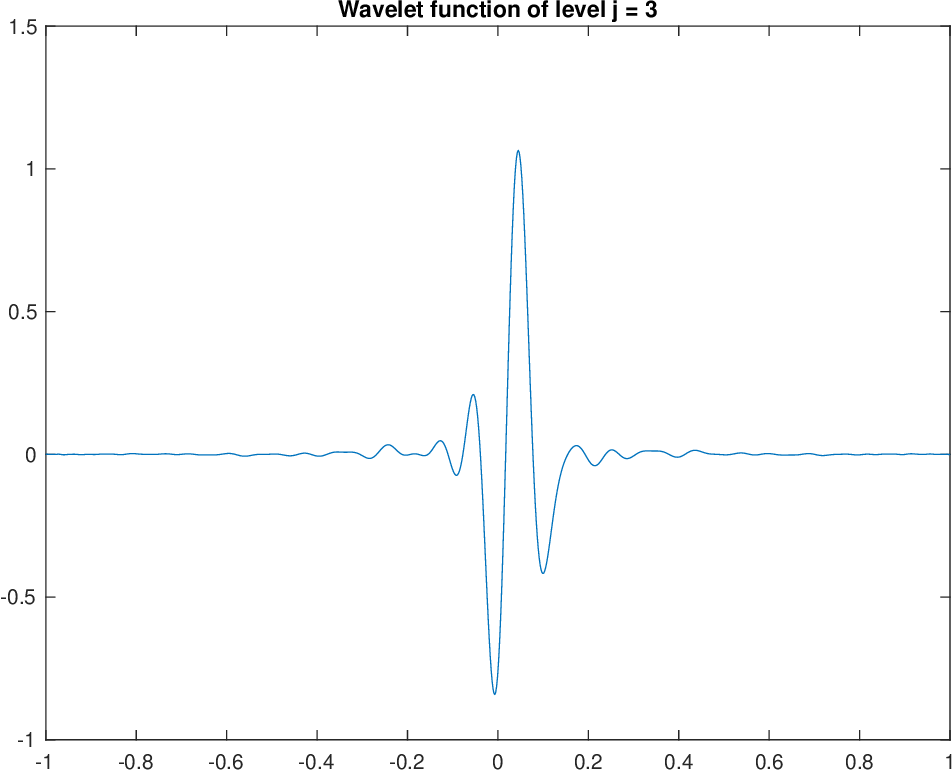}
\includegraphics[scale = 0.45]{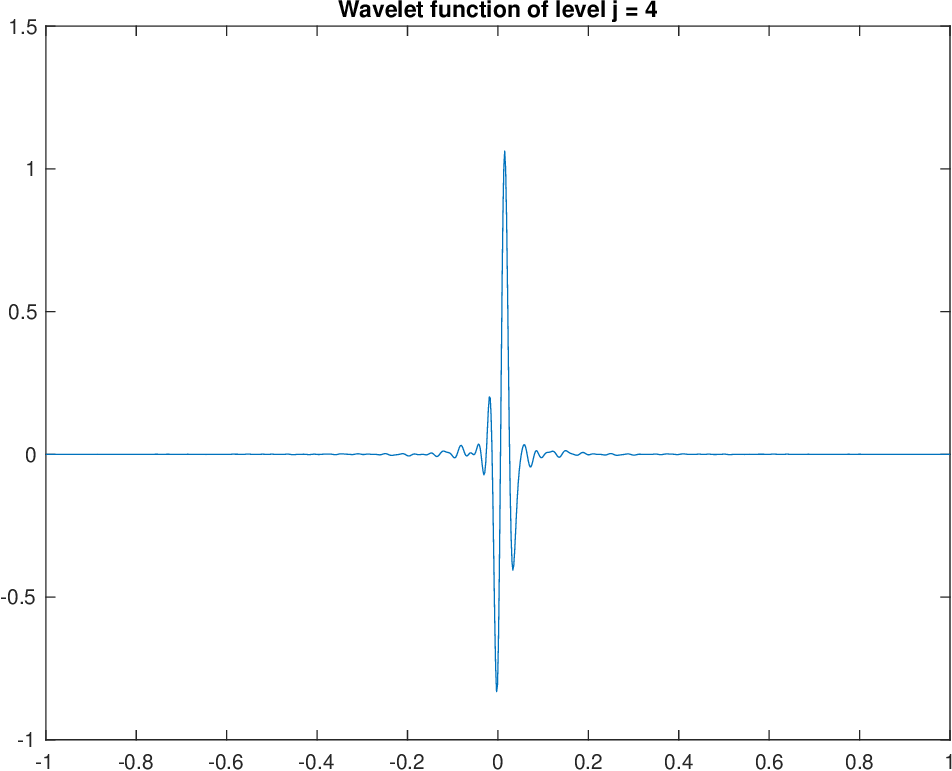}
\end{center}
\caption{ Wavelet functions for different values of $n = 3^j$ .\label{fig_wavelet_function_n} }
\end{figure}

In Figure \ref{fig_wavelet_functions_k}, we see each interpolating wavelet function $\psi_{n,k}^m$ is well localized around the Chebyshev node $y_k^{n}$.
Figure \ref{fig_wavelet_function_n} shows how this localization increases as $n$ increases. 

Below, we give other properties of the previous interpolating wavelets.
\begin{theorem}\label{th-wav1}
Let $n,m\in\NN$ be arbitrarily fixed with $n>m$. For all $k=1,\ldots,2n$, the first $n-m$ moments of the interpolating wavelets $\psi_{n,k}^m$ are zero, i.e.
\begin{equation}\label{mom-wav}
\int_{-1}^1 x^s\psi_{n,k}^m(x)w(x)dx=0\qquad s=0,1,\ldots, n-m, \qquad k=1,\ldots, 2n.
\end{equation}
Moreover, they form a Riesz basis of the details space $\W_n^m$ satisfying the following inequalities for  any vector $\vec{b}_{2n}=(b_{1}, \ldots, b_{2n})\in\RR^{2n}$ 
\begin{equation}
\label{Riesz-p-wav}
\C_1 \|\vec{b}_{2n}\|_{\ell^p}\le  \left\|\sum_{k=1}^{2n} b_{k}\psi_{n,k}^m\right\|_p\le \C(n,m)\ \|\vec{b}_{2n}\|_{\ell^p}, \qquad 1\le p\le\infty,
\end{equation}
where $\C(n,m)=\C_2(1+\log^2\frac nm)$, and $\C_1,\C_2>0$ are absolute constants  independent of $n, m, \vec{b}_{2n}$.
\end{theorem}
In the case that $n,m\in\NN$ satisfy \eqref{LCnm}, by \eqref{Riesz-p-wav} we get that 
\begin{equation}\label{norm-wav}
\|\psi_{n,k}^m\|_\infty\sim 1 \quad\mbox{and}\quad
\|\psi_{n,k}^m\|_p^p\sim \frac\pi{n}, \qquad 1\le p<\infty,
\end{equation}
hold for any $k=1,\ldots,2n$, where the notation $a\sim b$ means that there exist two absolute constants $C_1,C_2>0$ (independent of $n,m,k$) such that
$C_1 a\le b\le C_2 a$.

We conclude the section with the following theorem that states the analogous of \eqref{q-basis}--\eqref{qbasis-trans}, providing an orthogonal basis of the detail space $\W_n^m$ with the related change of basis formula.
\begin{theorem}\label{th-wort}
For all $n,m\in\NN$ with $m<n$, the polynomials 
\begin{equation}\label{wav-ort}
\psi_{n,r}^{m^{\hspace{.05cm}\mbox{\large $\bot$}}}(x)=\left\{\begin{array}
{ll}
\displaystyle \mu_{n, r}^m p_{2n-r}(x)+ \mu_{n, 2n-r}^m p_r(x) & n\le r<n+m\\ [.1in]
p_r(x)& n+m\le r\le 3n-m\\ [.1in]
\displaystyle \Phi_{3n,r}^\m(x)= \mu_{3n, r}^m p_{r}(x)- \mu_{3n, 6n-r}^m p_{6n-r}(x) & 3n-m< r<3n
\end{array}\right.
\end{equation}
constitute an orthogonal basis of 
$\W_n^m$ satisfying 
\[
<\psi_{n,r}^\m, \psi_{n,s}^\m>_{L_w^2}=v_{n,r}^m \delta_{r,s},\qquad r=n,\ldots, 3n-1 ,
\]
where
\begin{equation}\label{vnr}
 v_{n,r}^m=<\psi_{n,r}^\m
 \ \psi_{n,r}^\m >_{L^2_w}=\left\{\begin{array}
{ll}
\displaystyle \frac{m^2+(n-r)^2}{2m^2} & n< r<n+m\\ [.1in]
1& r=n \quad\mbox {or}\quad n+m\le r\le 3n-m\\ [.1in]
\displaystyle \frac{m^2+(3n-r)^2}{2m^2} & 3n-m< r<3n.
\end{array}\right. 
\end{equation}
In such a basis,  the interpolating wavelets  generating $\W_n^m$ can be expanded as follows
\begin{equation}\label{bchange-wav1}
\psi_{n,k}^m(x)=\sum_{r=n}^{3n-1}\rho_{r,k} \ \psi_{n,r}^{m^{\hspace{.05cm}\mbox{\large $\bot$}}}(x),
\qquad k=1,\ldots, 2n
\end{equation}
where
\begin{equation}\label{bchange-wav2}
\rho_{r,k}=\frac\pi{3n}\left\{\begin{array}{ll}
p_n(y_k^n) & r=n\\ [.1in]
p_r(y_k^n)+p_{|2n-r|}(y_k^n) & n<r\le 3n-m \quad r\ne 2n\\ [.1in]
p_{2n}(y_k^n)+\sqrt{2}p_0(y_k^n) & r=2n\\ [.1in]
p_r(y_k^n)+\mu_{n,r-2n}^m p_{r-2n}(y_k^n)-\mu_{n,4n-r}^m p_{4n-r}(y_k^n) & 3n-m<r<3n .
\end{array}\right.
\end{equation}
\end{theorem}
\section{Decomposition and reconstruction algorithms}
Due to \eqref{v+w}, for arbitrarily fixed $n,m\in\NN$ with $m<n$, any function $f_{3n}\in\V_{3n}^m$ can be uniquely decomposed as
\begin{equation}\label{f+g}
f_{3n}=f_n+g_{2n}, \qquad \mbox{with $\quad f_n\in\V_n^m\quad$ and 
$\quad g_{2n}\in\W_n^m$}.
\end{equation}
On the other hand, by using the scaling and wavelet bases introduced in the previous sections, the functions in \eqref{f+g} are uniquely determined by their respective scaling and wavelet coefficients.

In the following subsections we provide the decomposition and reconstruction formulas that lead to the fast computing of these coefficients. 

We start with the case that the orthogonal bases \eqref{q-basis} and \eqref{wav-ort} are fixed in $\V_n^m$ and $\W_n^m$ resp. Then we focus on the less trivial case where interpolating scaling and wavelet bases are used. 

In the former case, as expected, the formulas will be based on quasi orthogonal matrices. In the latter case, the formulas stated in \cite[Section 4]{CT-wave} will be generalized to our setting. Even if they do not involve orthogonal matrices, fast algorithms will be obtained also in this case.

\subsection{The orthogonal case}
With respect to the orthogonal bases \eqref{q-basis} and \eqref{wav-ort}  any function $f_n\in\V_n^m$ and any $g_{2n}\in\W_n^m$ are uniquely determined by the coefficient vectors 
\[
\vec{a}_n^\bot=(a_{0}^\bot,\ldots, a_{n-1}^\bot)\in\RR^n \qquad\mbox{and}\qquad\vec{b}_{2n}^\bot=(b_{n}^\bot,\ldots, b_{3n-1}^\bot)
\] 
defined by
\begin{eqnarray*}
f_n(x)&=&\sum_{k=0}^{n-1} a_{k}^\bot\Phi_{n,k}^\m(x)=
\vec{a}_n^\bot\cdot\vec{\Phi}_n^\bot,\qquad\qquad 
\vec{\Phi}_n^\bot=(\Phi_{n,1}^\m,\ldots, \Phi_{n,n}^\m)^T\\
g_{2n}(x)&=&\sum_{k=n}^{3n-1} b_{k}^\bot\psi_{n,k}^\m(x)
= \vec{b}_{2n}^\bot\cdot \vec{\psi}_n^\bot,\qquad\qquad
\vec{\psi}_n^\bot=(\psi_{n,n}^\m,\ldots,\psi_{n,3n-1}^\m)^T .
\end{eqnarray*}
The next theorem provides the decomposition and reconstruction formulas relating the vector $(\vec{a}_n^\bot, \vec{b}_{2n}^\bot)$ to the vector $\vec{a}_{3n}^\bot$ defining $f_{3n}=f_n+g_{2n}$ w.r.t. the orthogonal basis of $\V_{3n}^n$. The matrices in such formulas are based on the following elements
\begin{eqnarray*}
s_{k,j}&=&<\Phi_{n,k}^\m, \ \Phi_{3n,j}^\m>_{L^2_w},
\qquad k=0,\ldots, n-1\qquad j=0,\ldots, 3n-1,\\
w_{k,j}&=& <\psi_{n,k}^\m, \ \Phi_{3n,j}^\m>_{L^2_w},
\qquad k=n,\ldots, 3n-1\qquad j=0,\ldots, 3n-1.
\end{eqnarray*}
By \eqref{q-basis} and \eqref{wav-ort} we easily get
\begin{equation}\label{skj}
s_{k,j}=\left\{
\begin{array}{lll}
\delta_{k,j} & 0\le k\le n-m & \qquad 0\le j\le 3n-m\\ [.1in]
\mu_{n,k}^m\delta_{j,k}-\mu_{n,2n-k}^m\delta_{j,2n-k}
& n-m<k<n+m& \qquad  0\le j\le 3n-m\\ [.1in]
0 & 0\le k\le n-1& \qquad 3n-m< j< 3n\
\end{array}\right.
\end{equation}
and 
\begin{equation}\label{wkj}
w_{k,j} =\left\{
\begin{array}{lll}
\mu_{n,k}^m\delta_{j,k}+\mu_{n,2n-k}^m\delta_{j,2n-k}
& n\le k<n+m& \qquad  0\le j\le 3n-m\\ [.1in]
\delta_{k,j} & n+m\le k<3n & \qquad 0\le j< 3n.
\end{array}\right.
\end{equation}
\begin{theorem}\label{th-ort}
Under the previous setting, the decomposition and reconstruction formulas are 
\begin{equation}\label{eq-ort}
(\vec{a}_n^\bot, \ \vec{b}_{2n}^\bot)=
\vec{a}_{3n}^\bot M^{-1}\qquad\mbox{and}\qquad 
\vec{a}_{3n}^\bot= (\vec{a}_n^\bot, \ \vec{b}_{2n}^\bot) M
\end{equation}
where the matrix $M\in\RR^{3n\times 3n}$ and its inverse $M^{-1}$ have the following entries at the row $(k+1)$ and the column $(j+1)$ for $ k,j=0,\ldots,3n-1$
\begin{equation}\label{M-ort}
M_{k,j}=\frac 1{\nu_{3n,j}} \left\{\begin{array}{ll}
s_{k,j} & \mbox{if $\quad 0\le k< n$}\\ [.1in]
w_{k,j} & \mbox{if $\quad n\le k< 3n$}
\end{array}
\right.\qquad\mbox{and}\qquad
(M^{-1})_{k,j}=\left\{\begin{array}{ll}
\displaystyle
\frac{s_{j,k}}{\nu_{n,j}} & \mbox{if $\quad 0\le j< n$}\\ [.2in]
\displaystyle \frac{w_{j,k}}{v_{n,j}}
 & \mbox{if $\quad n\le j< 3n$}
\end{array}
\right.
\end{equation}
with $\nu_{n,j}=\|\Phi_{n,j}^\m\|_2^2$ and $v_{n,j}=\|\psi_{n,j}^\m\|_2^2$ explicitly given in \eqref{q-prod} and \eqref{wav-ort}, respectively.
\end{theorem}
\subsection{The interpolating case}
Now let us consider the interpolating scaling and wavelet bases functions, i.e.
\[
\V_n^m=span\{\Phi_{n,k}^m : \ k=1,\ldots, n\}\qquad\mbox{and}\qquad
\W_n^m=span\{\psi_{n,k}^m : \ k=1,\ldots, 2n\} .
\]
In such bases any function $f_n\in\V_n^m$ and any $g_{2n}\in\W_n^m$ can be expanded as follows
\begin{eqnarray}\label{fn}
f_n(x)&=&\sum_{k=1}^n a_{n,k}\Phi_{n,k}^m(x),\qquad a_{n,k}=f_n(x_k^n),\\
\label{gn}
g_{2n}(x)&=&\sum_{k=1}^{2n} b_{n,k}\psi_{n,k}^m(x),\qquad b_{n,k}=g_{2n}(y_k^n),
\end{eqnarray}
and they are uniquely determined by the scaling and wavelet coefficient vectors: 
\[
\vec{a}_n=(a_{n,1},\ldots, a_{n,n})\in\RR^n,\qquad 
\vec{b}_{2n}=(b_{n,1},\ldots, b_{n,2n})\in\RR^{2n}
\]
As regards $f_{3n}=f_n+g_{2n}$, of course, \eqref{fn} holds with $n$ replaced by $3n$. However, using the decomposition \eqref{XY}, we decompose the scaling vector $\vec{a}_{3n}\in\RR^{3n}$, determining $f_{3n}$, into the vectors 
\[
\vec{a}^{\ \prime}_{3n}=(a'_{1},\ldots, a'_{n})\in\RR^n \quad\mbox{and}\qquad 
\vec{a}^{\ \prime\prime}_{3n}=(a^{\prime\prime}_{1},\ldots, a^{\prime\prime}_{2n})\in\RR^{2n}
\]
defined by
\begin{equation}\label{f3n}
f_{3n}(x)=\sum_{k=1}^n a'_{k}\Phi_{3n}^m(x_k^n,x)+
\sum_{k=1}^{2n}a^{\prime\prime}_{k}\Phi_{3n}^m(y_k^n,x).
\end{equation}
We remark that
\begin{equation}\label{XY1}
x_k^n=x_{3k-1}^{3n}, \qquad k=1,\ldots,n
\end{equation}
implies that
\begin{equation}\label{XY2}
a'_k=f_{3n}(x_k^n)=f_{3n}(x_{3k-1}^{3n})= a_{3n, 3k-1}, \qquad k=1,\ldots,n,
\end{equation}
so that the 2nd, 5th, ...,$(3n-1)$--th entry of $\vec{a}_{3n}$ form the vector $\vec{a}_{3n}^{\ \prime}$, while the remaining entries of $\vec{a}_{3n}$ form the vector $\vec{a}_{3n}^{\ \prime\prime}$.

The following theorem generalizes a previous result in \cite{CT-wave} and provides the decomposition and reconstruction formulas to compute the vectors $\vec{a}_n$ and $\vec{b}_{2n}$ from $\vec{a}_{3n}$ decomposed as in \eqref{f3n}, and vice versa.
\begin{theorem}\label{th-decrec}
For all pairs of positive integers $m<n$, the basis coefficients defined by \eqref{fn}--\eqref{f3n}, for $f_{3n}=f_n+g_{2n}$, are related as follows:

{\bf Reconstruction formulas:} With $\Phi_{n,r}^\m$ defined in \eqref{q-basis}, we have\newline 
$\bullet$ For $k=1,\ldots,n$
\begin{equation}\label{rec-a1}
a'_{k}=a_{n,k}-\frac\pi{n}
\sum_{s=1}^{2n} b_{n,s} \sum_{r=0}^{n-1}
p_r(x_k^n)\Phi_{n,r}^\m(y_s^n) .
\end{equation}
$\bullet$ For $k=1,\ldots,2n$
\begin{equation}\label{rec-a2}
a^{\prime\prime}_{k}=b_{n,k}+\frac\pi{n}
\sum_{s=1}^n a_{n,s}\sum_{r=0}^{n-1}
\Phi_{n,r}^\m(y_k^n)p_r(x_s^n).
\end{equation}
{\bf Decomposition formulas:} With $\nu_{n,r}^m$ defined in \eqref{q-prod}, we have\newline
$\bullet$ For $k=1,\ldots,n$
\begin{equation}\label{dec-a}
a_{n,k}=\frac\pi{3n}\left[
\sum_{s=1}^n a'_{s} \sum_{r=0}^{n-1}
\frac{p_r(x_k^n)p_r(x_s^n)}{\nu_{n,r}^m}+
\sum_{s=1}^{2n} a^{\prime\prime}_{s} \sum_{r=0}^{n-1}\frac{p_r(x_k^n)\Phi_{n,r}^\m(y_s^n)}{\nu_{n,r}^m}\right] .
\end{equation}
$\bullet$ For $k=1,\ldots,2n$
\begin{equation}\label{dec-b}
b_{n,k}=a^{\prime\prime}_{k}-\frac\pi{3n}\left[
\sum_{s=1}^n a'_{s}\sum_{r=0}^{n-1}
\frac{\Phi_{n,r}^\m(y_k^n)p_r(x_s^n)}{\nu_{n,r}^m}+
\sum_{s=1}^{2n} a^{\prime\prime}_{s}\sum_{r=0}^{n-1}\frac{\Phi_{n,r}^\m(y_k^n)\Phi_{n,r}^\m(y_s^n)}{\nu_{n,r}^m}\right] .
\end{equation}
\end{theorem}

Following the derivation described in \cite{CT-wave}, if we change the order of the summation in the previous 
formulas (\ref{dec-a})-(\ref{rec-a2}) 
then both the decomposition and reconstruction algorithms can be described in the following 4 steps where, for brevity, we omit the notation of the parameters $n,m$ that do not vary and set 
\[
a_k=a_{n,k},\qquad b_k=b_{n,k},\qquad
x_j:= x_j^n,\qquad y_j:=y_j^n,\qquad \nu_r:= \nu_{n,r}^m, \qquad q_r:= \Phi_{n,r}^\m.
\]
\noindent{\sc Decomposition
Algorithm:}\vspace{.3cm}\newline\framebox{\parbox{10.5cm}{$
\begin{array}{ll}
Step \ 1. &  \mbox{Compute}\
\displaystyle\alpha_r=\sum_{s=1}^{n}\frac{\pi}{3n}
a'_{s}p_r(x_{s}), \hspace{.8cm} r=0,\ldots,n-1\\
Step \ 2. &  \mbox{Compute}\
\displaystyle\beta_r=\sum_{s=1}^{2n} \frac{\pi}{3n}
a''_{s}q_r(y_{s}), \hspace{.7cm} r=0,\ldots,n-1\\
Step \ 3. &\mbox{Compute}\ \displaystyle
a_{k}=\sum_{r=0}^{n-1}\frac{\alpha_r+\beta_r}{\nu_{r}}
p_r(x_{k}), \hspace{.3cm} k=1,\ldots,n \\
Step \ 4. & \mbox{Compute}\ \displaystyle
b_{k}=\sum_{r=0}^{n-1}\frac{\alpha_r+\beta_r}{\nu_{r}}
q_r(y_{k}), \hspace{.4cm} k=1,\ldots,2n
\end{array}
$ }}\vspace{.8cm}\newline{\sc Reconstruction
Algorithm:}\vspace{.3cm}\newline\framebox{\parbox{10.5cm}{$
\begin{array}{ll}
Step \ 1. &  \mbox{Compute}\
\displaystyle\alpha_r=\sum_{s=1}^{n}\frac{\pi}{n}
a_{s}p_r(x_{s}), \hspace{.8cm} r=0,\ldots,n-1\\
Step \ 2. &  \mbox{Compute}\
\displaystyle\beta_r=\sum_{s=1}^{2n}\frac{\pi}{3n}
b_{s}q_r(y_{s}), \hspace{.3cm} r=0,\ldots,n-1\\
Step \ 3. &\mbox{Compute}\ \displaystyle
a'_{k}=a_{k}-3
\sum_{r=0}^{n-1}\beta_r
p_r(x_{k}), \hspace{.3cm} k=1,.. ,n \\
Step \ 4. & \mbox{Compute}\ \displaystyle
a''_{k}=b_{k}+\sum_{r=0}^{n-1}\alpha_r q_r(y_{k}),
\hspace{.4cm} k=1,\ldots,2n
\end{array}
$ }}\vspace{.7cm}
\newline  Because $p_r$ is the Chebyshev polynomial of degree $r$
and also $q_r=\Phi_{n,r}^\m$ is defined in (\ref{q-basis}) by means of Chebyshev polynomials, using the trigonometric form of the Chebyshev polynomials (\ref{cheb-pol}), all the steps of the previous decomposition and
reconstruction algorithm can be computed efficiently, i.e., using 
${\cal O}(n\log n)$ flops. \newline
Indeed, if we use the Matlab definition for the several kinds of DCT (discrete cosine transformation) then, in both the algorithms, Step 1 can be computed as a DCT of type 2 and length $n$ while Step 2 as a DCT of type 2 and length $2n$,
Step 3 can be based on the DCT of type 3 and length n, and, finally, Step 4 can be computed using a DCT of type 3 and length $2n$.

Finally, we point out that given the sampling at the Chebyshev nodes of order $N=n_0 3^{Jmax}$ with $n_0\in \NN$, the previous decomposition and reconstruction formulas can be iterated $Jmax$ times. Moreover, taking into account that the nodes do not depend on $m$, at each decomposition step we can also choose a different $m$, but we have to remember the choice we make since in reconstruction we have to use the same $m$. In particular, once the value of $\theta$ is fixed, $\theta\in ]0,1[$, at each resolution level $n$ we can choose $m=\lfloor \theta n\rfloor$.

\begin{figure}[!htb]%
\begin{center}
\includegraphics[scale = 0.45]{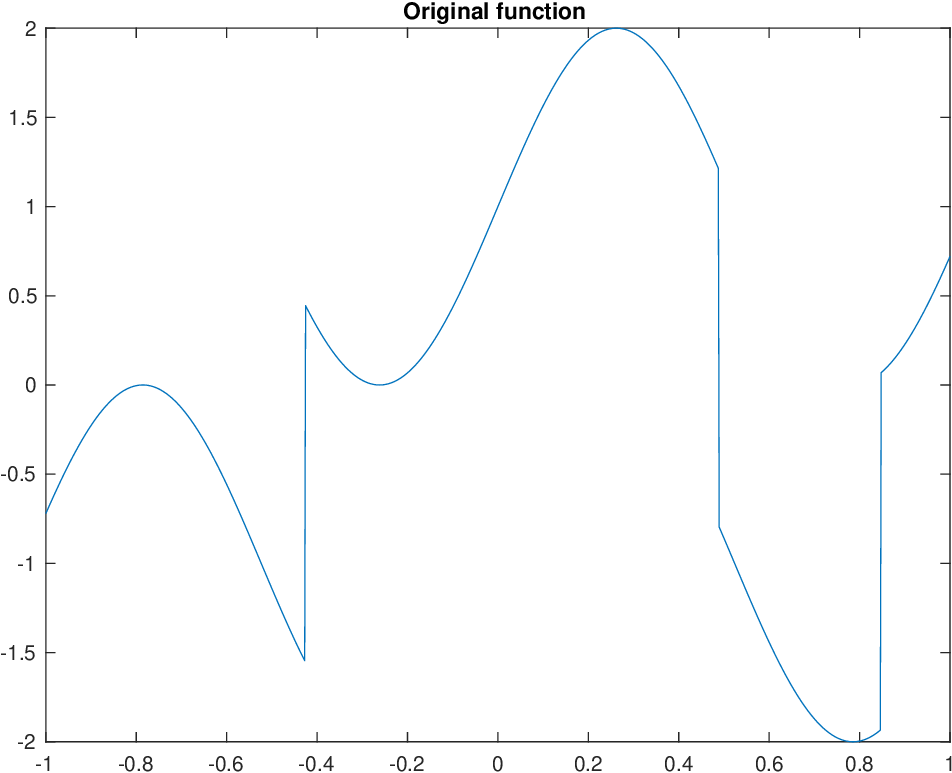}
\includegraphics[scale = 0.45]{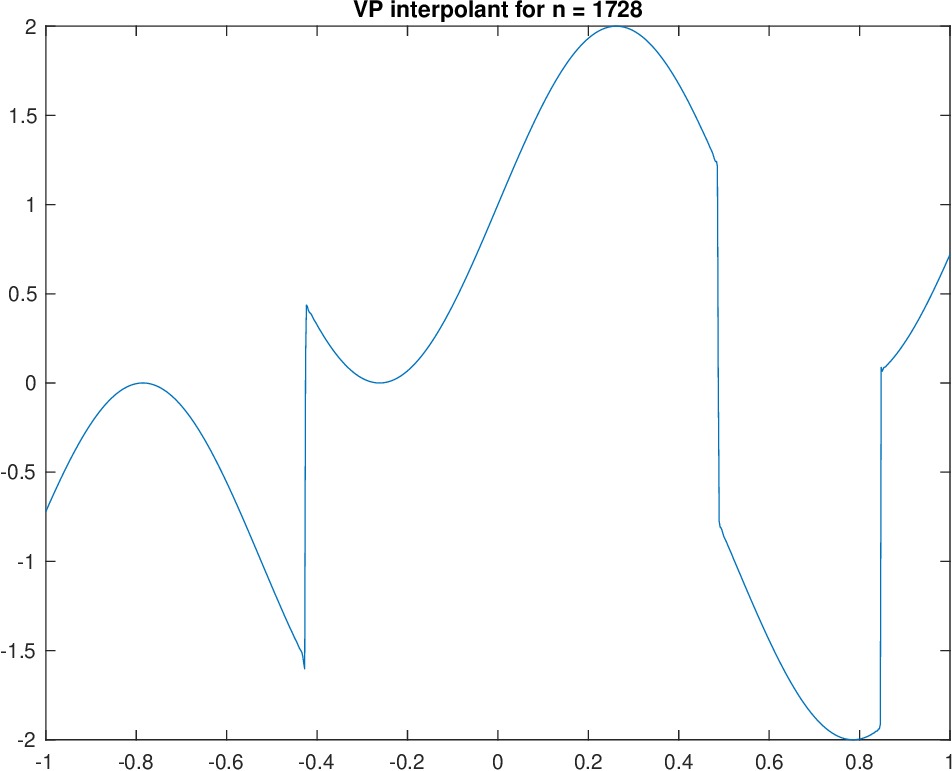}\\
\includegraphics[scale = 0.45]{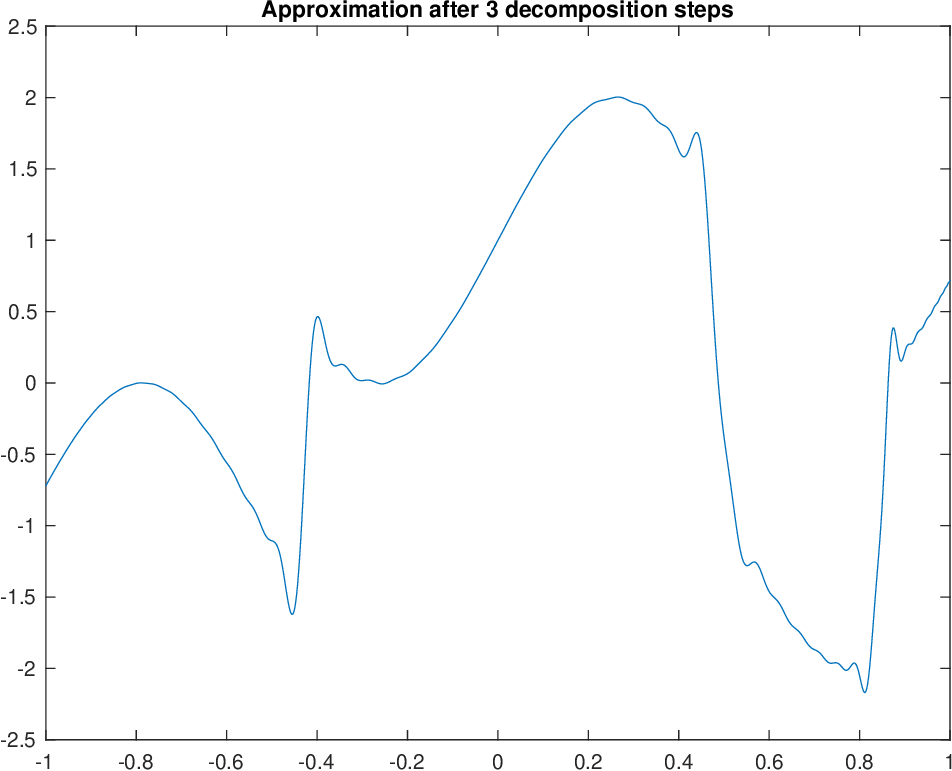}
\includegraphics[scale = 0.45]{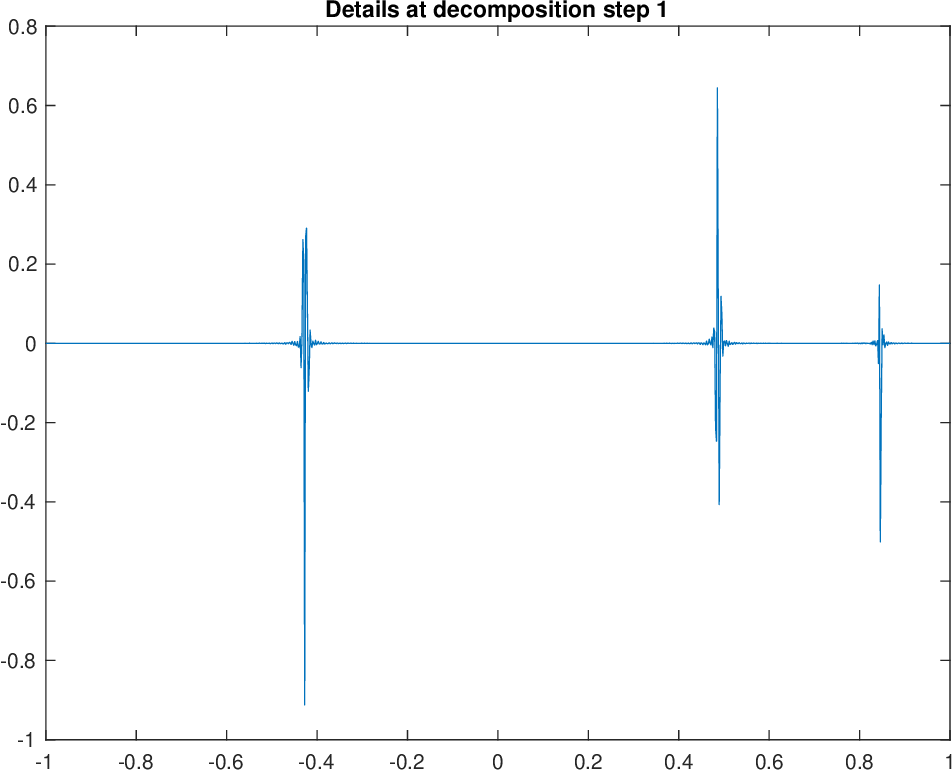}\\
\includegraphics[scale = 0.45]{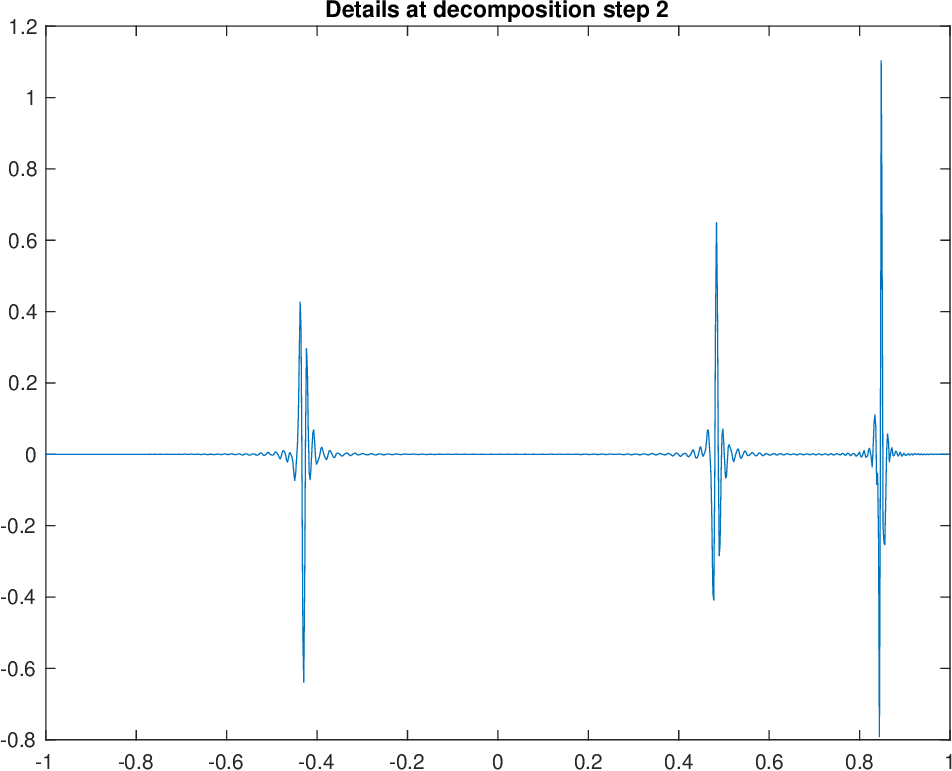}
\includegraphics[scale = 0.45]{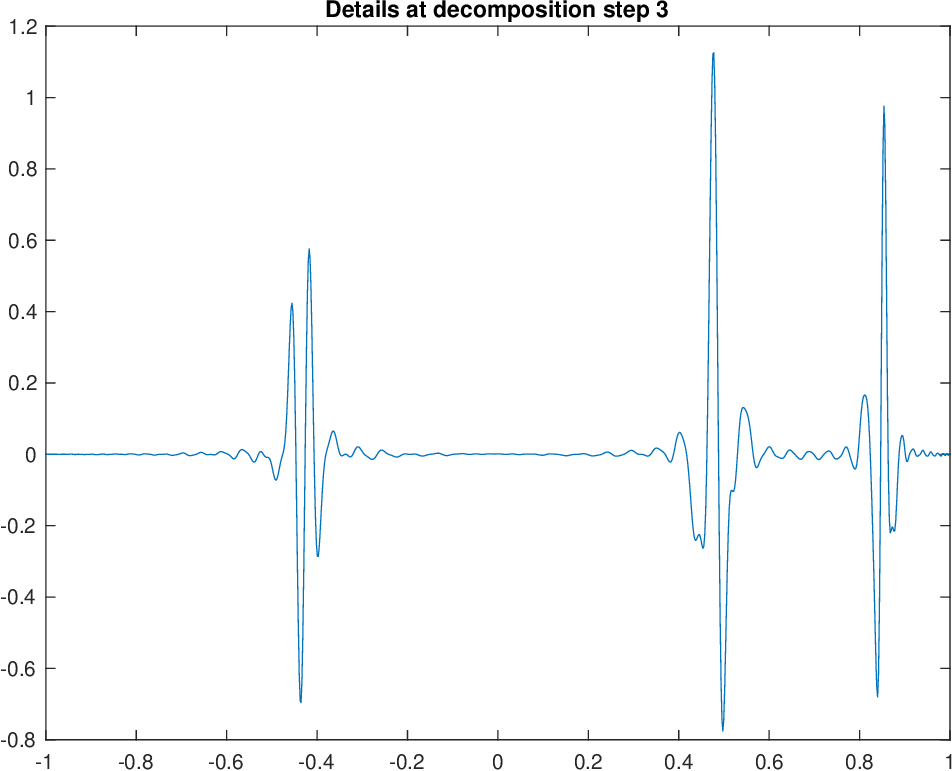}
\end{center}
\caption{ VP interpolation and wavelet decomposition for the function
$f(x) = \sin(6x)+\textrm{sign}(\sin(x+\exp(2x)))$
.\label{fig_wavelet_decomposition} }
\end{figure}

In Figure \ref{fig_wavelet_decomposition}, we illustrate an example of three decomposition steps for the function $f(x) = \sin(6x)+\textrm{sign}(\sin(x+\exp(2x)))$ by taking $m=\lfloor 0.7 n\rfloor$ at each level $n$.  
Starting from $n = 64 \cdot 3^3 = 1728$, we first project $f$ onto $\V_n^m$ by computing the VP interpolant $f_{1728}:=V_n^mf$ from the sampling of $f$ at $X_n$. Hence we iteratively decompose $f_{1728}$ as follows
\begin{eqnarray*}
 f_{1728}&=& f_{576}+g_{1152}\\
 &=& f_{192}+g_{384}+ g_{1152}\\
 &=& f_{64}+g_{128}+g_{384}+g_{1152}
\end{eqnarray*}
Note that we kept $\theta=0.7$
constant over the different levels but we could have chosen a different value for $\theta$ at each different level.

\section{Proofs}\label{Proofs}
{\large\bf Proof of Theorem \ref{th-sca}}\newline
$\bullet$ The identity $(a)$ trivially follows from \eqref{sca} and \eqref{VP-op}--\eqref{inva-sigma}.\newline
$\bullet$ To prove $(b)$, first we use \eqref{inva-sigma} and then apply the Gaussian rule
\begin{equation}\label{Gauss}
\int_{-1}^1 P(y)w(y)dy=\frac\pi n\sum_{k=1}^n P(x_k^n), \qquad \forall P\in\PP_{2n-1},
\end{equation}
to the polynomial $P(y)=y^s v_n^m(x,y)\in \PP_{s+n+m-1}\subseteq \PP_{2n-1}$, for all $x\in [-1,1]$, getting
\[
x^s=\int_{-1}^1 y^s v_n^m(x,y)w(y)dy=\frac\pi n\sum_{k=1}^n (x_k^n)^s v_n^m(x_k^n, x)=
\sum_{k=1}^n (x_k^n)^s \Phi_{n,k}^m(x).
\]
$\bullet$ To prove $(c)$, we note that $1=|\Phi_{n,k}^m(x_k^n)|\le \|\Phi_{n,k}^m\|_\infty$. On the other hand, by \eqref{VP-trig} we get the following trigonometric form of the scaling functions \cite{OT-APNUM21}
\begin{equation}\label{sca-trig}
\Phi_{n,.k}^m(x)=\frac 1{4mn}\left[\frac{\sin[m(t-t_k^n)]\sin[n(t-t_k^n)]}
{\sin^2[(t-t_k^n)/2]} + \frac{\sin[m(t+t_k^n)]\sin[n(t+t_k^n)]}
{\sin^2[(t+t_k^n)/2]}\right]
\end{equation}
where we set $t=\arccos x$ and $t_k^n=\arccos (x_k^n)$.
\newline
Hence, $(c)$ follows from \eqref{sca-trig} by taking into account that
\[
\left|\frac{\sin[N\theta]}{2N\sin[\theta/2]}\right|\le 1, \qquad \forall N\in\NN,\quad  \forall \theta\in (0,2\pi).
\]
$\bullet$ Regarding $(d)$, the case $p=1$ can be proved as follows by applying  (a) with $s=0$, and \eqref{c_inf}
\[
\frac\pi n=\int_{-1}^1\Phi_{n,k}^m(x)w(x)dx\le\int_{-1}^1|\Phi_{n,k}^m(x)|w(x)dx=\frac\pi n \int_{-1}^1|v_{n}^m(x_k^n, x)|w(x)dx\le \frac\pi n c_\infty.
\]
For $1<p<\infty$, the statement can be deduced from the (already proved) case $p=1$. More precisely, setting $\frac 1{p'}+\frac 1p=1$,  by H\"older's inequality we have
\[
\frac\pi n\le\int_{-1}^1|\Phi_{n,k}^m(x)|w(x)dx\le \left(\int_{-1}^1|\Phi_{n,k}^m(x)|^pw(x)dx\right)^\frac 1p \left(\int_{-1}^1w(x)dx\right)^\frac 1{p'}=
\|\Phi_{n,k}^m\|_p\  \pi^\frac 1{p^\prime}
\]
that proves the left-hand side inequality. Moreover, recalling $(c)$, we get
\[
\|\Phi_{n,k}^m\|_p^p=\int_{-1}^1|\Phi_{n,k}^m(x)|^pw(x)dx\le \int_{-1}^1|\Phi_{n,k}^m(x)|w(x)dx\le \frac\pi n c_\infty
\]
which concludes the proof of $(d)$.

$\bullet$ In order to prove $(e)$, we recall that for all algebraic polynomials $Q$ the following Marcinkiewicz--type inequality holds (see, e.g. \cite[Lemma 6.2]{OT-DRNA21})
\begin{equation}\label{Marci-n}
\left[\frac \pi n\sum_{k=1}^n |Q(x_k^n)|^p\right]^\frac 1p\le \left(1+2\pi \frac{\deg (Q)}n \right)\left[\int_{-1}^1 |Q(x)|^pw(x)dx\right]^\frac 1p, \qquad 1\le p<\infty.
\end{equation}
By \eqref{Marci-n} and \eqref{c_inf} we get $(e)$ as follows
\begin{eqnarray*}
\sum_{k=1}^n|\Phi_{n,k}^m(x)|&=&\frac \pi n \sum_{k=1}^n|v_n^m(x_k^n, x)|\le \left[1+2\pi\ \frac{n+m-1}{n}\right]\int_{-1}^1 |v_n^m(y,x)|w(y)dy\\
&\le&  \left[1+2\pi\ \frac{n+m-1}{n}\right] c_\infty.
\end{eqnarray*}
$\bullet$ Finally, let us prove $(f)$. To this aim, we set for brevity
\[
F(x)= \sum_{k=1}^{n} a_{k}\Phi_{n,k}^m(x), \qquad x\in [-1,1],
\]
and observe that $a_{k}=F(x_k^n)$ holds for $k=1,...,n$ by virtue of \eqref{sca-int}          .

Consequently, by using $(e)$, we get \eqref{Riesz-inf} as follows 
\begin{eqnarray*}
 \|\vec{a}_n\|_{\ell^\infty}&=&\max_{1\le k\le n}|F(x_k^n)|\le  \max_{|x|\le 1}|F(x)|\le \max_{|x|\le 1}\sum_{k=1}^{n}| a_{k}\Phi_{n,k}^m(x)|\\
&\le&
 \|\vec{a}_n\|_{\ell^\infty}\left( \max_{|x|\le 1}\sum_{k=1}^{n}|\Phi_{n,k}^m(x)|\right)\le
  \|\vec{a}_n\|_{\ell^\infty} \left[1+2\pi\frac{n+m-1}{n}\right]c_\infty
\end{eqnarray*}
Moreover,  by applying \eqref{Marci-n} with $Q=F$, we have
\[
\|\vec{a}_n\|_{\ell^p}=\left(\frac\pi n\sum_{k=1}^n|F(x_k^n)|^p\right)^\frac 1p\le  \left(1+2\pi \frac{n+m-1}n \right)\|F\|_p\le (1+4\pi )\|F\|_p, \qquad 1\le p<\infty
\]
which proves the first inequality in \eqref{Riesz-p}.  To prove the second ones, we recall that the dual space of $ L^{p}_w$ is $ L^{p'}_w$ where $\frac 1p+\frac 1{p'}=1$. Consequently,  there exists a function $G\in L^{p'}_w$ such that
\[
\|F\|_p=\left(\int_{-1}^1 |F(x)|^pw(x)dx\right)^\frac 1p=\int_{-1}^1 F(x) G(x) w(x)dx, \qquad \|G\|_{p'}=1.
\]
Using the definitions of $F$ and $\sigma_n^m G$, this means that
\[
\|F\|_p= \frac \pi n\sum_{k=1}^n a_{k}\int_{-1}^1 v_n^m(x_k^n, x) G(x) w(x)dx=  \frac \pi n\sum_{k=1}^n a_{k}[\sigma_n^mG(x_k^n)].
\]
Hence, for $p=1$, we conclude
\[
\|F\|_1\le \max_k |\sigma_n^mG(x_k^n)|\left(\frac \pi n\sum_{k=1}^n |a_{k}|\right)\le \|\sigma_n^mG\|_\infty \|\vec{a}_n\|_{\ell^1}\le c_\infty\|\vec{a}_n\|_{\ell^1}.
\]
Similarly, for $1<p<\infty$, by H\"older's inequality and \eqref{Marci-n}, we get
\[
\|F\|_p\le\left(\frac \pi n\sum_{k=1}^n |a_{k}|^p\right)^\frac 1p\left(\frac\pi n \sum_{k=1}^n |\sigma_n^mG(x_k^n)|^{p'}\right)^\frac 1{p'}\le
\|\vec{a}_n\|_{\ell^p} (1+4\pi)\|\sigma_n^mG\|_{p'}\le \|\vec{a}_n\|_{\ell^p} (1+4\pi)c_{p'}.
\]
\Proofend

{\large\bf Proof of Theorem \ref{th-wav}}\newline
Since $\psi_{n,k}^m\in \V_{3n}^m$, by \eqref{V-equi} we get
\begin{eqnarray*}
\psi_{n,k}^m(x)&=& \sum_{j=1}^{3n}\psi_{n,k}^m(x_{j}^{3n})\Phi_{3n,j}^m(x)
=\sum_{j=1}^{3n}\psi_{n,k}^m(x_{j}^{3n})\Phi_{3n}^m(x_j^{3n}, x)\\
&=& \sum_{h=1}^{2n}\psi_{n,k}^m(y_{h}^{n})\Phi_{3n}^m(y_h^n, x)
+ \sum_{h=1}^{n}\psi_{n,k}^m(x_{h}^{n})\Phi_{3n}^m(x_h^n, x)
\end{eqnarray*}
and using \eqref{wav-int} we obtain
\begin{equation}\label{eq1}
\psi_{n,k}^m(x)=\Phi_{3n}^m(y_k^n, x)
+ \sum_{h=1}^{n}\psi_{n,k}^m(x_{h}^{n})\Phi_{3n}^m(x_h^n, x).
\end{equation}
Hence, \eqref{wav} follows once proved that
\begin{equation}\label{rel-int}
\psi_{n,k}^m(x_h^n)=-\Phi_{n,h}^m(y_k^n), \qquad h=1,\ldots, n, \qquad k=1,\ldots, 2n.
\end{equation}
This is a consequence of $\psi_{n,k}^m\bot\Phi_{n,h}^m$. Indeed,  using \eqref{eq1}, \eqref{inva-sigma}, and \eqref{sca-int} we have
\begin{eqnarray*}
0&=& < \Phi_{n,h}^m, \psi_{n,k}^m>_{L^2_w}\\
&=&
< \Phi_{n,h}^m,\ \Phi_{3n}^m(y_k^n, \cdot) >_{L^2_w}
+ \sum_{j=1}^{n}\psi_{n,k}^m(x_{j}^{n}) < \Phi_{n,h}^m,\ \Phi_{3n}^m(x_j^n, \cdot) >_{L^2_w}\\
&=&\frac\pi{3n}\left[\sigma_{3n}^m[\Phi_{n,h}^m](y_k^n)+
\sum_{j=1}^{n}\psi_{n,k}^m(x_{j}^{n}) \sigma_{3n}^m[\Phi_{n,h}^m](x_j^n)
\right]\\
&=&\frac\pi{3n}\left[\Phi_{n,h}^m(y_k^n)+
\sum_{j=1}^{n}\psi_{n,k}^m(x_{j}^{n}) \Phi_{n,h}^m(x_j^n)
\right]\\
&=&\frac\pi{3n}\left[\Phi_{n,h}^m(y_k^n)+
\psi_{n,k}^m(x_{h}^{n})\right]
\end{eqnarray*}
that yields \eqref{rel-int}.\Proofend

{\large\bf Proof of Theorem \ref{th-wav1}}\newline
The identity \eqref{mom-wav} trivially follows from  $\psi_{n,k}^m\bot \V_n^m$ and $\PP_{n-m}\subset\V_n^m$.

In order to prove \eqref{Riesz-p-wav}, we note that the polynomial
\[
F(x)=\sum_{k=1}^{2n}b_{k}\psi_{n,k}^m(x), \qquad x\in [-1,1]
\]
belongs to $\V_{3n}^m$. Hence, by \eqref{V-equi}, it can be written as
\[
F(x)=V_{3n}^mF(x)=\sum_{j=1}^{3n}F(x_j^{3n})\Phi_{3n,j}^m(x), \qquad x\in [-1,1]
\]
where, due to \eqref{wav-int} and \eqref{rel-int}, we have
\begin{equation}\label{Fj}
F(x_j^{3n})=\left\{
\begin{array}{lll}
F(y_h^n)=b_{h} & \mbox{if $\ x_j^{3n}=y_h^n$}, & h=1,\ldots, 2n\\ [.1in]
\displaystyle F(x_h^n)= - \sum_{k=1}^{2n}b_{k}\Phi_{n,h}^m(y_k^n) & \mbox{if $\ x_j^{3n}=x_h^n$},& h=1,\ldots, n.
\end{array}
\right.
\end{equation}
On the other hand, if we apply \eqref{Marci-n} replacing $n$ with $3n$, and use $Y_{2n}\subset X_{3n}$, then for any algebraic polynomial $Q$ we get
\begin{equation}\label{Marci-n-wav}
\left[\frac \pi {3n}\sum_{k=1}^{2n} |Q(y_k^n)|^p\right]^\frac 1p\le \left(1+2\pi \frac{\deg (Q)}{3n} \right)\left[\int_{-1}^1 |Q(x)|^pw(x)dx\right]^\frac 1p, \qquad 1\le p<\infty.
\end{equation}
In particular, \eqref{Marci-n-wav} implies that
\begin{eqnarray*}
\sum_{k=1}^{2n}|\Phi_{n,h}^m(y_k^n)|&=&\frac \pi{n}
\sum_{k=1}^{2n}|v_n^m(x_h^n, y_k^n)|\le 3\left(1+2\pi \frac{n+m-1}{3n} \right) \int_{-1}^1 |v_n^m(x_h^n,x)|w(x)dx\\
&\le& (3+4\pi)c_\infty(n,m)
\end{eqnarray*}
and 
\[
|F(x_h^n)|=\left|\sum_{k=1}^{2n}b_{k}\Phi_{n,h}^m(y_k^n)\right|\le \max_{1\le k\le 2n}|b_k|
\left(\sum_{k=1}^{2n}|\Phi_{n,h}^m(y_k^n)|\right)\le (3+4\pi)c_\infty(n,m)\|\vec{b}_{2n}\|_{\ell^\infty}
\]
hold for any $h=1,\ldots, n$.

Using these estimates and Thm. \ref{th-sca} we get 
\begin{eqnarray*}
\|\vec{b}_{2n}\|_{\ell^\infty}&=&\max_{1\le h\le 2n}|F(y_h^n)|\le \|F\|_\infty= \max_{|x|\le 1}\left|
\sum_{j=1}^{3n}F(x_j^{3n})\Phi_{3n,j}^m(x)\right|\\
&\le& (1+4\pi) c_\infty(3n,m)\left(\max_{1\le j\le 3n}|F(x_j^{3n})|\right)\\
&=&  (1+4\pi) c_\infty(3n,m)\max\left\{\max_{1\le h\le 2n}|F(y_h^{n})|, \ \max_{1\le h\le n}|F(x_h^n)|\right\}\\
&\le& (1+4\pi) c_\infty(3n,m)\max\left\{\|\vec{b}_{2n}\|_{\ell^\infty}, \ (3+4\pi)c_\infty(n,m)\|\vec{b}_{2n}\|_{\ell^\infty}\right\}\\
&\le& (3+4\pi)^2 c_\infty(3n,m)c_\infty(n,m)\|\vec{b}_{2n}\|_{\ell^\infty},
\end{eqnarray*}
that proves \eqref{Riesz-p-wav} for $p=\infty$.

Now, let us prove \eqref{Riesz-p-wav} when $1\le p<\infty$. For brevity, in the sequel we denote by $\C$ a positive constant that can take different values at different occurrences but it is always independent of $n,m, \vec{b}_{2n}$.

The left-hand side inequality in \eqref{Riesz-p-wav}
can be deduced from \eqref{Fj}, $Y_{2n}\subset X_{3n}$, and the left-hand side inequality in \eqref{Riesz-p}, as follows
\[
\|\vec{b}_{2n}\|_{\ell^p}=\left(\frac\pi{2n}\sum_{h=1}^{2n}|F(y_h^n)|^p\right)^\frac 1p\le 
\\C\left(\frac\pi{3n}\sum_{j=1}^{3n}|F(x_j^{3n})|^p\right)^\frac 1p\le
\C\|F\|_p.
\]
Finally, let us prove the right-hand side inequality in \eqref{Riesz-p-wav} for any $1\le p<\infty$. 

Note that the right-hand side inequality in \eqref{Riesz-p} and \eqref{Fj} imply
\begin{eqnarray*}
\|F\|_p&=& \left\|\sum_{j=1}^{3n}F(x_j^{3n})\Phi_{3n,j}^m\right\|_p\\
&\le& (1+4\pi)\ c_{p'}(3n,m) 
\left(\frac\pi{3n}\sum_{j=1}^{3n}|F(x_j^{3n})|^p\right)^\frac 1p\\
&\le& \C \ c_{p'}(3n,m)  \left[\left(\frac\pi{3n}\sum_{h=1}^{2n}|F(y_h^{n})|^p\right)^\frac 1p
+\left(\frac\pi{3n}\sum_{h=1}^{n}|F(x_h^{n})|^{p}\right)^\frac 1{p}\right]\\
&\le&\C \ c_{p'}(3n,m)  \left[\left(\frac\pi{2n}\sum_{h=1}^{2n}|b_{h}|^p\right)^\frac 1p
+\left(\frac\pi{n}\sum_{h=1}^{n}|F(x_h^{n})|^{p}\right)^\frac 1{p}\right].
\end{eqnarray*}
Hence, to complete the proof we are going to prove that
\[
\left(\frac\pi n\sum_{h=1}^{n}|F(x_h^{n})|^{p}\right)^\frac 1{p}\le 
\C\ c_{p'}(n,m) \left(\frac\pi {2n} \sum_{k=1}^{2n}|b_{k,n}|^p\right)^\frac 1p.
\]
To this aim, first we note that \eqref{Fj} and \eqref{Marci-n} imply
\begin{eqnarray*}
\left(\frac\pi n\sum_{h=1}^n|F(x_h^n)|^p\right)^\frac 1p&=&\left(\frac \pi n\sum_{h=1}^n\left|\sum_{k=1}^{2n}b_{k}\Phi_{n}^m(x_h^n, y_k^n)\right|^p\right)^\frac 1p\\
&\le& \C \left(\int_{-1}^1 \left|\sum_{k=1}^{2n}b_{k}\Phi_{n}^m(x, y_k^n)\right|^pw(x)dx\right)^\frac 1p.
\end{eqnarray*}
Then, setting $\frac 1p +\frac 1{p'}=1$, and recalling that
\[
\forall f\in L^p_w, \quad \exists g\in L^{p'}_w\quad s.t. \qquad \|g\|_{p'}=1\quad\mbox{and}\quad \|f\|_p=\int_{-1}^1f(x)g(x)w(x)dx
\]
by using \eqref{VP-op}, H\"older's inequality,  \eqref{Marci-n-wav} and \eqref{LCp}, we conclude
\begin{eqnarray*}
\left(\frac\pi n\sum_{h=1}^n|F(x_h^n)|^p\right)^\frac 1p
&\le& \C \left(\int_{-1}^1 \left|\sum_{k=1}^{2n}b_{k}\Phi_{n}^m(x, y_k^n)\right|^pw(x)dx\right)^\frac 1p\\
&=&\C \int_{-1}^1 \left(\sum_{k=1}^{2n}b_{k}\Phi_{n}^m(x, y_k^n)\right)g(x)w(x)dx\\
&=& \C \left[\frac\pi n \sum_{k=1}^{2n}b_{k}\ 
\sigma_{n}^m g (y_k^n)\right]\\
&\le& \C\left(\frac\pi n \sum_{k=1}^{2n}|b_{k}|^p\right)^\frac 1p\left(\frac\pi n \sum_{k=1}^{2n}\left|\sigma_{n}^m g (y_k^n)\right|^{p'}\right)^\frac 1{p'}\\
&\le& \C\left(\frac\pi {2n} \sum_{k=1}^{2n}|b_{k}|^p\right)^\frac 1p\ 
\left\|\sigma_{n}^m g \right\|_{p'}\\
&\le& \C\left(\frac\pi {2n} \sum_{k=1}^{2n}|b_{k}|^p\right)^\frac 1p
\ c_{p'}(n,m) .
\end{eqnarray*}
\Proofend 

{\large \bf Proof of Theorem \ref{th-wort}}

By the definition, for any $r,s=n,\ldots, 3n-1$, we trivially deduce that
\[
<\psi_{n,r}^{m^{\hspace{.05cm}\mbox{\large $\bot$}}},
 \ \psi_{n,s}^{m^{\hspace{.05cm}\mbox{\large $\bot$}}}>_{L^2_w}=\delta_{r,s}\cdot \left\{\begin{array}
{ll}
\displaystyle \frac{m^2+(n-r)^2}{2m^2} & n< r<n+m\\ [.1in]
1& r=n \quad\mbox {or}\quad n+m\le r\le 3n-m\\ [.1in]
\displaystyle \frac{m^2+(3n-r)^2}{2m^2} & 3n-m< r<3n.
\end{array}\right.
\]
Moreover, for any $r=n,\ldots, 3n-1$, the polynomial 
$\psi_{n,r}^{m^{\hspace{.05cm}\mbox{\large $\bot$}}}$ belongs to 
$\W_n^m$ because, first of all, it can be trivially expanded in terms of the orthogonal polynomials spanning $\V_{3n}^m$
\[
\Phi_{3n,s}^{m^{\hspace{.05cm}\mbox{\large $\bot$}}}(x)=\left\{\begin{array}{ll}
p_s(x) & 0\le s\le 3n-m\\ [.1in]
\displaystyle \mu_{3n, s}^m p_{s}(x)- \mu_{3n, 6n-s}^m p_{6n-s}(x) & 3n-m< s<3n
\end{array}\right.
\]
and hence $\psi_{n,r}^{m^{\hspace{.05cm}\mbox{\large $\bot$}}}\in \V_{3n}^m$. Furthermore, $\psi_{n,r}^\m$ is orthogonal to $\V_n^m$ since we have
\[
<\psi_{n,r}^{m^{\hspace{.05cm}\mbox{\large $\bot$}}}, \ 
\Phi_{n,s}^{m^{\hspace{.05cm}\mbox{\large $\bot$}}}>_{L^2_w}=0, \qquad r=n,\ldots, 3n-1,\qquad s=0,\ldots, n-1
\]
Taking into account that $\dim \W_n^m=2n$, the previous arguments prove the first part of the statement, namely $\{\psi_{n,r}^{m^{\hspace{.05cm}\mbox{\large $\bot$}}}, \ r=n,\ldots, 3n-1\}$ is an orthogonal basis of $\W_n^m$.

Now, let us prove \eqref{bchange-wav1}. 

For any  $k=n,\ldots, 3n-1$, taking into account that $\psi_{n,k}^m\in \V_{3n}^m$,  by \eqref{Vnm-q}, \eqref{wav-int}, \eqref{rel-int}, and \eqref{Vnm}, we have
\begin{eqnarray*}
\psi_{n,k}^m(x)&=&\sum_{r=0}^{3n-1}\Phi_{3n,r}^{\m}(x)\left[\frac\pi{3n}
\sum_{h=1}^{3n}\psi_{n,k}^m(x_h^{3n})p_r(x_h^{3n})\right]\\
&=&\frac\pi{3n}\sum_{r=0}^{3n-1}\Phi_{3n,r}^\m(x)\left[
\sum_{h=1}^{2n}\psi_{n,k}^m(y_h^{n})p_r(y_h^{n})
+
\sum_{h=1}^{n}\psi_{n,k}^m(x_h^{n})p_r(x_h^{n})\right]\\
&=&\frac\pi{3n}\sum_{r=0}^{3n-1}\Phi_{3n,r}^\m(x)\left[
p_r(y_k^{n})-\sum_{h=1}^{n}\Phi_{n,h}^m(y_k^{n})p_r(x_h^{n})\right]\\
&=&\frac\pi{3n}\sum_{r=0}^{3n-1}\Phi_{3n,r}^\m(x)\left[
p_r(y_k^{n})-V_n^mp_r(y_k^{n})\right].
\end{eqnarray*}
Hence, by  \eqref{inva} and \eqref{Vnm-q}, we obtain
\begin{eqnarray*}
\psi_{n,k}^m(x)&=&\frac\pi{3n}\sum_{r=n-m+1}^{3n-1}\Phi_{3n,r}^\m(x)\left[p_r(y_k^n)-V_n^mp_r(y_k^n)\right]\\
&=&\frac\pi{3n}\sum_{r=n-m+1}^{3n-1}\Phi_{3n,r}^\m(x)
\left[p_r(y_k^n)-\sum_{s=0}^{n-1}\Phi_{n,s}^\m(y_k^n)\Sigma_{r,s}\right]
\end{eqnarray*}
where, for brevity, we set
\[
\Sigma_{r,s}=\frac\pi n\sum_{h=1}^np_r(x_h^n)p_s(x_h^n),\qquad r,s\in\NN.
\]
Taking into account that $p_n(x_h^n)=0$, $p_{2n}(x_h^n)=-\sqrt{2}p_0(x_h^n)$, and that
\begin{eqnarray*}
p_{2n-s}(x_h^n)&=&-p_{s}(x_h^n), \qquad s=1,\ldots,n-1\\
p_{2n+s}(x_h^n)&=&-p_{s}(x_h^n), \qquad s=1,\ldots,n-1
\end{eqnarray*}
hold for all $h=1,\ldots, n$, we can easily prove
\begin{equation}\label{sigma}
\Sigma_{r,s}=\left\{\begin{array}{clll}
<p_s,p_r>_{L^2_w} &=\delta_{s,r} & \mbox{if}& 0\le r<n\\
0 & & \mbox{if}&r=n\\
<p_s, -p_{2n-r}>_{L^2_w} &= - \delta_{s,2n-r} &\mbox{if}& n< r<2n\\
<p_s, -\sqrt{2}p_0>_{L^2_w} & =-\sqrt{2}\delta_{s,0} &\mbox{if}& r=2n\\
<p_s,-p_{r-2n}>_{L^2_w} &=-\delta_{s,r-2n} &\mbox{if}& 2n< r<3n
\end{array}\right.\qquad s=0,\ldots, n-1.
\end{equation}
Consequently, by using \eqref{sigma} and making explicit the form of the orthogonal scaling functions, simple computations yield the following identities 
\begin{eqnarray*}
S_1&:=&\sum_{r=n-m+1}^{n+m-1}\Phi_{3n,r}^\m(x)
\left[p_r(y_k^n)-\sum_{s=0}^{n-1}\Phi_{n,s}^\m(y_k^n)\Sigma_{r,s}\right]
\\
&=& \sum_{r=n-m+1}^{n+m-1}p_r(x)
\left[p_r(y_k^n)-\sum_{s=0}^{n-1}\Phi_{n,s}^\m(y_k^n)\Sigma_{r,s}\right]\\
&=& \sum_{r=n-m+1}^{n-1}p_r(x)\left[p_r(y_k^n)-\Phi_{n,r}^\m(y_k^n)\right]+
p_n(x)p_n(y_k^n)+\sum_{r=n+1}^{n+m-1}p_r(x)\left[p_r(y_k^n)+\Phi_{n,2n-r}^\m(y_k^n)\right]\\
&=&p_n(x)p_n(y_k^n)+ \left\{\sum_{r=n-m+1}^{n-1}+\sum_{r=n+1}^{n+m-1}\right\}
p_r(x)\mu_{n,2n-r}^m\left[p_r(y_k^n)+p_{2n-r}(y_k^n)\right]
\\
&=& p_n(x)p_n(y_k^n)\sum_{s=n+1}^{n+m-1}p_{2n-s}(x)\mu_{n,s}^m
\left[p_{2n-s}(y_k^n)+p_s(y_k^n)\right]+
\sum_{r=n+1}^{n+m-1}p_r(x)\mu_{n,2n-r}^m\left[p_r(y_k^n)+p_{2n-r}(y_k^n)\right]
\\
&=&p_n(x)p_n(y_k^n) +\sum_{r=n+1}^{n+m-1}\left[p_r(y_k^n)+p_{2n-r}(y_k^n)\right]
\left[ \mu_{n,r}^mp_{2n-r}(x)+\mu_{n,2n-r}^mp_r(x)\right].
\end{eqnarray*}
\begin{eqnarray*}
S_2&:=&\sum_{r=n+m}^{3n-m}\Phi_{3n,r}^\m(x)
\left[p_r(y_k^n)-\sum_{s=0}^{n-1}\Phi_{n,s}^\m(y_k^n)\Sigma_{r,s}\right]\\
&=& \sum_{r=n+m, \ r\ne 2n}^{3n-m}
p_r(x)\left[p_r(y_k^n)+\Phi_{n,|2n-r|}^\m(y_k^n)\right]+ p_{2n}(x)[p_{2n}(y_k^n)+\sqrt{2}p_0(y_k^n)]\\
&=& \sum_{r=n+m, \ r\ne 2n}^{3n-m}
p_r(x)\left[p_r(y_k^n)+p_{|2n-r|}^m(y_k^n)\right]+ p_{2n}(x)[p_{2n}(y_k^n)+\sqrt{2}p_0(y_k^n)].
\end{eqnarray*}
\begin{eqnarray*}
S_3&:=&\sum_{r=3n-m+1}^{3n-1}\Phi_{3n,r}^\m(x)
\left[p_r(y_k^n)-\sum_{s=0}^{n-1}\Phi_{n,s}^\m(y_k^n)\Sigma_{r,s}\right]\\
&=&\sum_{r=3n-m+1}^{3n-1}\Phi_{3n,r}^\m(x) \left[p_r(y_k^n)+\Phi_{n,r-2n}^\m(y_k^n)\right]\\
&=&\sum_{r=3n-m+1}^{3n-1}\left[\mu_{3n,r}^m p_r(x)- \mu_{3n,6n-r}^m p_{6n-r}(x)\right] \left[p_r(y_k^n)+\mu_{n,r-2n}^mp_{r-2n}(y_k^n)-
\mu_{n,4n-r}^mp_{4n-r}(y_k^n)\right].
\end{eqnarray*}
In conclusion, taking into account that
\[
\psi_{n,k}^m(x)=\frac\pi{3n}\left(S_1+S_2+S_3\right)
\]
the previous identities imply \eqref{bchange-wav1}--\eqref{bchange-wav2}. \Proofend

{\large\bf Proof of Theorem \ref{th-ort}}\newline
First of all note that, due to orthogonality, we have the following two--scale relations
\begin{eqnarray}
\Phi_{n,k}^\m(x)&=& \sum_{j=0}^{3n-1}\frac{<\Phi_{n,k}^\m, \ \Phi_{3n,j}^\m>_{L^2_w}}{\|\Phi_{3n,j}^\m\|_2^2}
\ \Phi_{3n,j}^\m(x), \qquad |x|\le 1\\
 \nonumber
&& \qquad\qquad\qquad\qquad \qquad\qquad k=1,\ldots,n\\
\psi_{n,k}^\m(x)&=& \sum_{j=0}^{3n-1}\frac{<\psi_{n,k}^\m, \ \Phi_{3n,j}^\m>_{L^2_w}}{\|\Phi_{3n,j}^\m\|_2^2}
 \Phi_{3n,j}^\m(x), \qquad |x|\le 1\\
 \nonumber
&&\qquad\qquad\qquad \qquad\qquad\qquad k=n,\ldots,3n-1 
\end{eqnarray}
that yield the following  transformations of the basis vectors
\begin{equation}\label{2scale-ort}
\left(\begin{array}{c}
\vec{\Phi}_n^\bot\\ [.1in]
\vec{\psi}_n^\bot
\end{array}\right) = M \vec{\Phi}_{3n}^\bot \qquad \mbox{and} \qquad
 \vec{\Phi}_{3n}^\bot = M^{-1} \left(\begin{array}{c}
\vec{\Phi}_n^\bot\\ [.1in]
\vec{\psi}_n^\bot
\end{array}\right)
\end{equation}
where we set
\[
\vec{\Phi}_n=\left(\begin{array}{c}
\Phi_{n,1}^m(x)\\
\vdots\\
\vdots\\
\Phi_{n,n}^m(x)
\end{array}\right)\qquad \vec{\psi}_n=\left(\begin{array}{c}
\psi_{n,1}^m(x)\\
\vdots\\
\vdots\\
\psi_{n,2n}^m(x)
\end{array}\right)
\]
Note that non singularity of $M$ and the explicit form of $M^{-1}$ can be easily deduced by introducing the following diagonal matrices
\[
\Delta=diag (\|\Phi_{3n,j}^\m\|_2)_j, \qquad D_1=diag (\|\Phi_{n,k}^\m\|_2 ) _k,\qquad D_2= diag (\|\psi_{n,k}^\m\|_2 )_k 
\]
and taking into account that the matrix 
\[
Q=D^{-1} M\Delta, \qquad D=\left(\begin{array}{c|c}
  D_1 &  0\\
  \hline
  0   &    D_2
\end{array}\right)
\]
is an orthogonal matrix in $\RR^{3n\times 3n}$, since its entries are 
\[
Q_{k,j}=\left\{\begin{array}{ll}
\displaystyle <\ \frac{\Phi_{n,k}^\m}{\|\Phi_{n,k}^\m\|_2}, \  \frac{\Phi_{3n,j}^\m}{\|\Phi_{3n,j}^\m\|_2}\ >_w,  &  k=0,\ldots,n-1 \\ [.4in]
\displaystyle <\ \frac{\psi_{n,k}^\m}{\|\psi_{n,k}^\m\|_2}, \  \frac{\Phi_{3n,j}^\m}{\|\Phi_{3n,j}^\m\|_2}\ >_w,    & k=n,\ldots, 3n-1
\end{array}\right.\quad j=0,\ldots, 3n-1
\]
and trivially satisfy
\[
\sum_{j=0}^{3n-1}Q_{k,j} Q_{h,j}=\delta_{h,k}, \qquad h,k=0,\ldots, 3n-1
\]
In conclusion, by applying \eqref{2scale-ort}, we can deduce the decomposition and reconstruction formulas \eqref{eq-ort} from 
\begin{eqnarray*}
f_n+g_{2n}&=&(\vec{a}_n^\bot, \vec{b}_{2n}^\bot)\cdot \left(\begin{array}{c}
\vec{\Phi}_n^\bot\\ [.1in]
\vec{\psi}_n^\bot
\end{array}\right)=(\vec{a}_n^\bot, \vec{b}_{2n}^\bot) M \vec{\Phi}_{3n}^\bot\\
f_{3n}&=& \vec{a}_{3n}^\bot \vec{\Phi}_{3n}^\bot= \vec{a}_{3n}^\bot M^{-1} \left(\begin{array}{c}
\vec{\Phi}_n^\bot\\ [.1in]
\vec{\psi}_n^\bot
\end{array}\right).
\end{eqnarray*}
by equating the coefficients of the respective basis functions.
\Proofend

{\large\bf Proof of Theorem \ref{th-decrec}}\newline
The proof can be achieved as in \cite{CT-wave}. For completeness, we report here the main steps adapted to our notation.

We start from the following two--scale relations 
\begin{eqnarray}
\label{2scale-1}
\Phi_{n,k}^m(x)&=&\Phi_{3n}^m(x_k^n,x) + \sum_{h=1}^{2n} \Phi_{n,k}^m(y_h^{n})\Phi_{3n}^m(y_h^n,x), \qquad k=1,\ldots,n,\\
\label{2scale-2}
\psi_{n,k}^m(x)&=&\Phi_{3n}^m(y_k^n,x) - \sum_{h=1}^n \Phi_{n,h}^m(y_k^{n})\Phi_{3n}^m(x_h^n,x), \qquad k=1,\ldots,2n
\end{eqnarray}
that in matrix form can be written as
\begin{equation}\label{scale-vect}
\left(\begin{array}{c}
\vec{\Phi}_n\\ [.2cm]
\vec{\Psi}_n
\end{array}\right)=\left(\begin{array}{c|c}
I & A_n\\ [.2cm] \hline
&\\[-.2cm]
-A_n^T & I
\end{array}\right)
\left(\begin{array}{c}
\vec{\Phi}'_{3n}\\ [.2cm]
\underline{\Phi}''_{3n}
\end{array}\right),
\end{equation}
where $A_n\in \RR^{n\times 2n}$ denotes the matrix with entries
\[
A_{h,k}=\Phi_{n,h}(y_k^n), \qquad h=1,\ldots,n,\qquad k=1,\ldots, 2n,
\]
$I$ denotes a rectangular identity matrix  of suitable dimensions with entries $I_{h,k}=\delta_{h,k}$, and we set
\[
\vec{\Phi}_n=\left(\begin{array}{c}
\Phi_{n,1}^m(x)\\
\vdots\\
\vdots\\
\Phi_{n,n}^m(x)
\end{array}\right)\quad \vec{\psi}_n=\left(\begin{array}{c}
\psi_{n,1}^m(x)\\
\vdots\\
\vdots\\
\psi_{n,2n}^m(x)
\end{array}\right)\quad
\vec{\Phi}_{3n}'=\left(\begin{array}{c}
\Phi_{3n}^m(x_1^n, x)\\
\vdots\\
\vdots\\
\Phi_{3n}^m(x_n^n, x)
\end{array}\right)\quad
\vec{\Phi}_{3n}''=\left(\begin{array}{c}
\Phi_{3n}^m(y_1^n, x)\\
\vdots\\
\vdots\\
\Phi_{3n}^m(y_{2n}^n, x)
\end{array}\right).
\]
From \eqref{scale-vect} we deduce
\[
f_n+g_{2n}=(\vec{a}_n\ \vec{b}_{2n}) 
\left(\begin{array}{c}
\vec{\Phi}_n\\ [.2cm]
\vec{\Psi}_n
\end{array}\right)=
(\vec{a}_n\ \vec{b}_{2n}) 
\left(\begin{array}{c|c}
I & A_n\\ [.2cm] \hline
&\\[-.2cm]
-A_n^T & I
\end{array}\right)
\left(\begin{array}{c}
\vec{\Phi}'_{3n}\\ [.2cm]
\underline{\Phi}''_{3n}
\end{array}\right),
\]
that compared with
\[
f_n+g_{2n}=f_{3n}= (\vec{a}_{3n}'\ \vec{a}_{3n}'' )
\left(\begin{array}{c}
\vec{\Phi}'_{3n}\\ [.2cm]
\underline{\Phi}''_{3n}
\end{array}\right),
\]
yields the reconstruction formulas \eqref{rec-a1}--\eqref{rec-a2}. 

Since the decomposition formulas \eqref{dec-a}--\eqref{dec-b} can be similarly deduced from the inverse two scale relations 
\[
\left(\begin{array}{c}
\vec{\Phi}^\prime_{3n}\\ [.2cm]
\underline{\Phi}^{\prime\prime}_{3n}
\end{array}\right)
=M^{-1}
\left(\begin{array}{c}
\vec{\Phi}_n\\ [.2cm]
\vec{\Psi}_n
\end{array}\right),\qquad \quad\mbox{where}\qquad M=\left(\begin{array}{c|c}
I & A_n\\ [.2cm] \hline
&\\[-.2cm]
-A_n^T & I
\end{array}\right)
\]
we need to investigate the form of $M^{-1}$.  To this aim, setting
\[
G_n=I+A_nA_n^T
\]
we can easily check that
\begin{equation}\label{M-1}
M^{-1}=\left(\begin{array}{c|c}
G_n^{-1} & -G_n^{-1}A_n\\ [.2cm] \hline
&\\[-.2cm]
A_n^TG_n^{-1} & I-A_n^TG_n^{-1}A_n
\end{array}\right).
\end{equation}
Looking for a more explicit form of $M^{-1}$, let us observe that,
 by \eqref{Gauss}, \eqref{XY}, and \eqref{sca-int}, we get
\begin{eqnarray*}
<\Phi_{n,r}^m, \Phi_{n,s}^m>_{L^2_w}&=&\frac\pi{3n}
\sum_{j=1}^{3n}\Phi_{n,r}^m(x_j^{3n})\Phi_{n,s}^m(x_j^{3n})\\
&=& \frac\pi{3n}\left[\sum_{k=1}^n\Phi_{n,r}^m(x_k^{n})\Phi_{n,s}^m(x_k^{n})+\sum_{k=1}^{2n}\Phi_{n,r}^m(y_k^{n})\Phi_{n,s}^m(y_k^{n})\right]\\
&=& \frac\pi{3n}\left[\delta_{r,s}+\sum_{k=1}^{2n}A_{r,k}A_{s,k}\right]
\end{eqnarray*}
which proves that the entries of the matrix $G_n$ are
the following
\[
(G_n)_{r,s}=(I+A_nA_n^T)_{r,s}=\frac{3n}\pi <\Phi_{n,r}^m, \Phi_{n,s}^m>_{L^2_w}, \qquad r,s=1,\ldots,n.
\]
Thus, by using \eqref{qbasis-trans}, we obtain
\[
(G_n)_{r,s}=\frac{3n}\pi <\Phi_{n,r}^m, \Phi_{n,s}^m>_{L^2_w}
=\frac{3n}\pi \left(\frac\pi n\right)^2 \sum_{i=0}^{n-1}p_i(x_r^n)p_i(x_s^n)<\Phi_{n,i}^\m, \Phi_{n,i}^\m>_{L^2_w} 
\]
i.e.
\begin{equation}\label{G}
(G_n)_{r,s}= 3 \left[\frac\pi n \sum_{i=0}^{n-1} \nu_{n,i}^m \ p_i(x_r^n)p_i(x_s^n)\right]
\qquad r,s=1,\ldots,n.
\end{equation}
Using this identity and taking into account that 
\[
\frac\pi n \sum_{k=1}^{n}p_r(x_k^n)p_s(x_k^n)=\delta_{r,s}=\frac\pi n \sum_{i=0}^{n-1}p_i(x_r^n)p_i(x_s^n), \qquad r,s=1,\ldots,n,
\]
we can easily check that the entries of $G_n^{-1}$ are
\begin{equation}\label{G-1}
(G_n^{-1})_{r,s}=\frac 13 \left[\frac\pi n\sum_{i=0}^{n-1} \frac 1{\nu_{n,i}^m} \ p_i(x_r^n)p_i(x_s^n)\right],
\qquad r,s=1,\ldots,n.
\end{equation}
Moreover, by \eqref{G-1} and \eqref{qbasis-trans}, we also get
that the other block matrices in \eqref{M-1} have the following entries
\begin{eqnarray}
\label{GA}
(G_n^{-1}A_n)_{r,s}&=&\frac\pi{3n}\sum_{i=0}^{n-1}
\frac 1{\nu_{n,i}^m} \ p_i(x_r^n)\Phi_{n,i}^\m (y_s^n),
\quad r=1,\ldots,n, \quad s=1,\ldots,2n\\
\label{AG}
(A_n^TG_n^{-1})_{r,s}&=&\frac\pi{3n}\sum_{i=0}^{n-1}
\frac 1{\nu_{n,i}^m} \ p_i(x_s^n)\Phi_{n,i}^\m (y_r^n),
\quad r=1,\ldots,2n, \quad s=1,\ldots,n\\
\label{AGA}
(A_n^TG_n^{-1}A_n)_{r,s}&=&\frac\pi{3n}\sum_{i=0}^{n-1}
\frac 1{\nu_{n,i}^m} \ \Phi_{n,i}^\m(y_r^n)\Phi_{n,i}^\m (y_s^n),
\qquad r,s=1,\ldots,2n.
\end{eqnarray}
In conclusion, by \eqref{M-1}, the inverse two--scale relations are
\[
\left(\begin{array}{c}
\vec{\Phi}^\prime_{3n}\\ [.2cm]
\underline{\Phi}_{3n}^{\prime\prime}
\end{array}\right)=
\left(\begin{array}{c|c}
G_n^{-1} & -G_n^{-1}A_n\\ [.2cm] \hline
&\\[-.2cm]
A_n^TG_n^{-1} & I-A_n^TG_n^{-1}A_n
\end{array}\right)
\left(\begin{array}{c}
\vec{\Phi}_n\\ [.2cm]
\vec{\Psi}_n
\end{array}\right)
\]
and imply
\[
(\vec{a}_n\ \vec{b}_{2n})= (\vec{a}_{3n}^\prime\ \vec{a}_{3n}^{\prime\prime})\left(\begin{array}{c|c}
G_n^{-1} & -G_n^{-1}A_n\\ [.2cm] \hline
&\\[-.2cm]
A_n^TG_n^{-1} & I-A_n^TG_n^{-1}A_n
\end{array}\right)
\]
which, by means of \eqref{G-1}--\eqref{AGA}, yields the decomposition formulas \eqref{dec-a}--\eqref{dec-b}.
\Proofend
\section{Conclusions}
Without any mother function to dilate and translate, we present a family of scaling and wavelet functions on $[-1,1]$ that are non-orthogonal, interpolating, and very smooth being polynomials. 

The scaling functions at a resolution level $n$ are defined as fundamental VP polynomials interpolating at Chebyshev zeros of the first kind and order $n$, whose good approximation properties are already known in the literature \cite{Th-2012}. They generate a non-standard multiresolution analysis where instead of duplicating, we triplicate the resolution degree $n$, introducing $2n$ wavelets spanning the detail space at the resolution level $n$. Such wavelets are also  polynomials interpolating at those Chebyshev nodes of the scaling functions at level $3n$ which are missing at the resolution level $n$. 

Both the scaling and wavelet functions have the peculiarity of depending on a free integer parameter $m$ that can be arbitrarily chosen in the set $\{1,\ldots, n-1\}$. In particular, this parameter determines the number of null moments for wavelets of level $n$ (the first $n-m$ moments). 
For the special choice $n=2\cdot 3^j$ and $m=3^j$, $j\in\NN$, we get the scaling and wavelet functions introduced in \cite{CT-wave} from which many ideas were taken.

Besides several theoretical results, Riesz stability is proved for scaling functions as well as for wavelets w.r.t. the uniform norm and any $L^p$ norm weighted by the Chebyshev weight. Moreover, fast decomposition and reconstruction algorithms, based on DCT, are given.

 Going along the same line, the proposed procedure could be extended to any odd scaling factor replacing the higher resolution level $3n$ with $5n$, $7n$, etc.... Moreover, similar results could be derived for the other Chebyshev weights of the second, third and fourth kind, as done in \cite{CT-wave} for a special parameter setting.

\bibliographystyle{abbrv}
\bibliography{biblioW}

\begin{thebibliography}{10}

\bibitem{CT-wave}
M.~Capobianco and W.~Themistoclakis.
\newblock Interpolating polynomial wavelets on [-1, 1].
\newblock {\em Advances in Computational Mathematics}, 23(4):353 – 374, 2004.

\bibitem{Chui-book}
C.~Chui.
\newblock {\em An Introduction to Wavelets}.
\newblock Academic Press, London, 1992.

\bibitem{Cohen1993}
A.~Cohen, I.~Daubechies, and P.~Vial.
\newblock Wavelets on the interval and fast wavelet transforms.
\newblock 1(1):54 – 81, 1993.

\bibitem{Dau-book}
I.~Daubechies.
\newblock {\em Ten Lectures on Wavelets}.
\newblock SIAM, Philadelphia, 1992.

\bibitem{BOT-Prandtl}
M.~De~Bonis, D.~Occorsio, and W.~Themistoclakis.
\newblock Filtered interpolation for solving {P}randtl’s integro-differential
  equations.
\newblock {\em Numerical Algorithms}, 88(2):679 – 709, 2021.

\bibitem{DT}
Z.~Ditzian and V.~Totik.
\newblock {\em Moduli of smoothness}.
\newblock Springer-Verlag, New York, 1987.

\bibitem{FisPre}
B.~Fischer and J.~Prestin.
\newblock Wavelets based on orthogonal polynomials.
\newblock 66(220):1593 – 1618, 1997.

\bibitem{FT-wave}
B.~Fischer and W.~Themistoclakis.
\newblock Orthogonal polynomial wavelets.
\newblock {\em Numerical Algorithms}, 30(1):37 – 58, 2002.

\bibitem{Gus}
J.~Gustavsson.
\newblock On interpolation of weighted $\ell^p$ spaces and {O}vchinnikov's
  theorem.
\newblock {\em Stud. Math.}, 72:237 – 251, 1982.

\bibitem{KilPre}
T.~Kilgore and J.~Prestin.
\newblock Polynomial wavelets on the interval.
\newblock 12(1):95 – 110, 1996.

\bibitem{MT-air}
G.~Mastroianni and W.~Themistoclakis.
\newblock A numerical method for the generalized airfoil equation based on the
  de la {V}all\'ee {P}oussin interpolation.
\newblock {\em J. Comput. Appl. Math.}, 180:71--105, 2005.

\bibitem{ORT-Etna}
D.~Occorsio, G.~Ramella, and W.~Themistoclakis.
\newblock Filtered polynomial interpolation for scaling 3{D} images.
\newblock {\em Electronic Transactions on Numerical Analysis}, 59:295 – 318,
  2023.

\bibitem{ORT-Jmiv}
D.~Occorsio, G.~Ramella, and W.~Themistoclakis.
\newblock Image scaling by de la {V}all\'ee {P}oussin filtered interpolation.
\newblock {\em Journal of Mathematical Imaging and Vision}, 65(3):513 – 541,
  2023.

\bibitem{ORT-RemSen}
D.~Occorsio, G.~Ramella, and W.~Themistoclakis.
\newblock An open image resizing framework for remote sensing applications and
  beyond.
\newblock {\em Remote Sensing}, 15(16), 2023.

\bibitem{ORT-Hilbert1}
D.~Occorsio, M.~Russo, and W.~Themistoclakis.
\newblock Filtered integration rules for finite weighted {H}ilbert transforms.
\newblock {\em Journal of Computational and Applied Mathematics}, 410(114166),
  2022.

\bibitem{ORT-Hilbert2}
D.~Occorsio, M.~G. Russo, and W.~Themistoclakis.
\newblock Filtered integration rules for finite weighted {H}ilbert transforms
  ii.
\newblock {\em Dolomites Research Notes on Approximation}, 15(3):93 – 104,
  2022.

\bibitem{ORT-Cauchy}
D.~Occorsio, M.~G. Russo, and W.~Themistoclakis.
\newblock On solving some {C}auchy singular integral equations by de la
  {V}all\'ee {P}oussin filtered approximation.
\newblock {\em Applied Numerical Mathematics}, 200:358 – 378, 2024.

\bibitem{OT-APNUM21}
D.~Occorsio and W.~Themistoclakis.
\newblock On the filtered polynomial interpolation at {C}hebyshev nodes.
\newblock {\em Appl. Numer. Math.}, 166:272--287, 2021.

\bibitem{OT-DRNA21}
D.~Occorsio and W.~Themistoclakis.
\newblock Some remarks on filtered polynomial interpolation at {C}hebyshev
  nodes.
\newblock {\em Dolomites Research Notes on Approximation}, 14(2):68--84, 2021.

\bibitem{PSTasche}
G.~Plonka, K.~Selig, and M.~Tasche.
\newblock On the construction of wavelets on a bounded interval.
\newblock 4(1):357 – 388, 1995.

\bibitem{Th-1999}
W.~Themistoclakis.
\newblock Some interpolating operators of de la {V}all\'ee {P}oussin type.
\newblock {\em Acta Mathematica Hungarica}, 84(3):221 – 235, 1999.

\bibitem{Th-2012}
W.~Themistoclakis.
\newblock Uniform approximation on [-1, 1] via discrete de la {V}all\'ee
  {P}oussin means.
\newblock {\em Numerical Algorithms}, 60(4):593 – 612, 2012.

\bibitem{Th-L1}
W.~Themistoclakis.
\newblock Weighted {L}1 approximation on [-1,1] via discrete de la {V}all\'ee
  {P}oussin means.
\newblock {\em Mathematics and Computers in Simulation}, 147:279 – 292, 2018.

\bibitem{TB-GVP}
W.~Themistoclakis and M.~Van~Barel.
\newblock Generalized de la {V}all\'ee {P}oussin approximations on [-1, 1].
\newblock {\em Numerical Algorithms}, 75(1):1 – 31, 2017.

\bibitem{Zyg-book}
A.~Zygmund.
\newblock {\em Trigonometric series}.
\newblock Cambridge University Press, 2015.

\end{thebibliography}

\section*{Funding declaration}
W. Themistoclakis was partially supported by INdAM-GNCS.
M. Van Barel was partially supported by the Fund for Scientific Research–Flanders (Belgium), project G0B0123N.

\end{document}